\documentclass[11pt]{article}
\usepackage{amsfonts}
\usepackage{color}
\usepackage{amsmath,amssymb,comment}
\newtheorem{theo}{Theorem}[section]
\newtheorem{lem}[theo]{Lemma}
\newtheorem{cor}[theo]{Corollary}
\newtheorem{prop}[theo]{Proposition}

\newcommand{\mysection}[1]{\section{#1} \setcounter{equation}{0}}
\newcommand{\proof}{{\sc Proof.} \quad}
\newcommand{\proofc}{{\sc Proof} \ }
\newcommand{\be}{\begin{equation} \label}
\newcommand{\ee}{\end{equation}}
\newcommand{\bea}{\begin{eqnarray}\label}
\newcommand{\eea}{\end{eqnarray}}
\newcommand{\bas}{\begin{eqnarray*}}
\newcommand{\eas}{\end{eqnarray*}}
\newcommand{\bit}{\begin{itemize}}
\newcommand{\eit}{\end{itemize}}
\newcommand{\qed}{\hfill$\Box$ \vskip.2cm}
\newcommand{\nn}{\nonumber}
\newcommand{\R}{\mathbb{R}}
\newcommand{\N}{\mathbb{N}}
\newcommand{\pO}{\partial\Omega}

\newcommand{\eps}{\varepsilon}

\newcommand{\supp}{{\rm supp} \, }

\newcommand{\wto}{\rightharpoonup}
\newcommand{\wsto}{\stackrel{\star}{\rightharpoonup}}
\newcommand{\hra}{\hookrightarrow}
\newcommand{\io}{\int_\Omega}
\newcommand{\na}{\nabla}
\newcommand{\Del}{\Delta}

\newcommand{\pa}{\partial}
\newcommand{\bom}{\overline{\Omega}}
\newcommand{\Om}{\Omega}

\newcommand{\ov}{\overline}
\newcommand{\wh}{\widehat}

\newcommand{\hs}{\hspace*}

\newcommand{\vp}{\varphi}
\newcommand{\lbal}{\left\{ \begin{array}{l}}
\newcommand{\lball}{\left\{ \begin{array}{ll}}
\newcommand{\ear}{\end{array} \right.}
\newcommand{\ouz}{\ov{u}_0}

\newcommand{\abs}{\\[5pt]}

\newcommand{\tme}{T_{max,\eps}}

\newcommand{\ueps}{u_\eps}
\newcommand{\veps}{v_\eps}
\newcommand{\weps}{w_\eps}

\newcommand{\yeps}{y_\eps}

\newcommand{\ig}{\int_G}
\begin{document}
\enlargethispage{10mm}
\title{A quantitative strong parabolic maximum principle\\
and application to a taxis-type migration-consumption model\\
involving signal-dependent degenerate diffusion}
\author{
Michael Winkler\footnote{michael.winkler@math.uni-paderborn.de}\\
{\small Universit\"at Paderborn, Institut f\"ur Mathematik,}\\
{\small 33098 Paderborn, Germany} }
\date{}
\maketitle
\begin{abstract}
\noindent 
  The taxis-type migration-consumption model accounting for signal-dependent motilities, as given by
  \bas
	\left\{ \begin{array}{l}
	u_t = \Delta \big(u\phi(v)\big), \\[1mm]
	v_t = \Delta v-uv,
	\ear
	\qquad \qquad (\star)
  \eas
  is considered for suitably smooth functions $\phi:[0,\infty)\to\R$ which are such that $\phi>0$ on $(0,\infty)$, but 
  that in addition $\phi(0)=0$ with $\phi'(0)>0$.\abs
  In order to appropriately cope with the diffusion degeneracies thereby included,
  this study separately examines
  the Neumann problem for the linear equation
  \bas
	V_t  = \Del V + \na\cdot \big( a(x,t)V\big) + b(x,t)V
  \eas
  and establishes a statement on how pointwise positive lower bounds for nonnegative solutions 
  depend on the supremum and the mass of the initial data, and on integrability features 
  of $a$ and $b$.\abs
  This is thereafter used as a key tool in the derivation of a result on global existence of solutions to ($\star$), 
  smooth and classical for positive times, under the mere assumption that the suitably regular initial data be nonnegative
  in both components.
  Apart from that, these solutions are seen to stabilize toward some equilibrium,	
  and as a qualitative effect genuinely due to degeneracy in diffusion, a criterion on initial smallness of the
  second component is identified as sufficient for this limit state to be spatially nonconstant.\abs
\noindent {\bf Key words:} maximum principle; functional inequality; chemotaxis; pattern formation\\
 {\bf MSC 2020:} 35K55 (primary); 35B50, 35B40, 35K10, 92C17 (secondary)
\end{abstract}
\newpage
\section{Introduction}\label{intro}
The primary subject of this study is the initial-boundary value problem
\be{0}
	\left\{ \begin{array}{ll}
	u_t = \Delta \big(u\phi(v)\big),
	\qquad & x\in\Om, \ t>0, \\[1mm]
	v_t = \Delta v-uv,
	\qquad & x\in\Om, \ t>0, \\[1mm]
	\frac{\partial u}{\partial\nu}=\frac{\partial v}{\partial\nu}=0,
	\qquad & x\in\pO, \ t>0, \\[1mm]
	u(x,0)=u_0(x), \quad v(x,0)=v_0(x),
	\qquad & x\in\Om,
	\end{array} \right.
\ee
in a smoothly bounded domain $\Om\subset\R^n$, $n\ge 1$, with nonnegative initial data $u_0$ and $v_0$,
and with a suitably regular nonnegative function $\phi$ on $[0,\infty)$.
Parabolic systems of this form arise in the modeling of collective behavior in bacterial ensembles, as represented through
their population densities $u=u(x,t)$, under the influence of certain signal substances, described by their concentrations
$v=v(x,t)$.
Here the approach to describe migration by second-order operators of the form in the first equation from (\ref{0})
accounts for recent advances in the modeling literature addressing bacterial movement especially in situations in which 
cell motility {\em as a whole} may be biased by chemical cues (\cite{fu}, \cite{liu}). \abs
Due to the resulting precise quantitative connection between diffusive and cross-diffusive contributions, (\ref{0}) can
be viewed as singling out a special subclass of Keller-Segel type chemotaxis models with their commonly much less restricted interdependencies (\cite{hillen_painter2009}).
The ambition to understand implications of the particular structural properties going along with this type of link,
and, in more general, to carve out possible peculiarities and characteristic features in population models including
migration operators of this form, has stimulated considerable activity in the mathematical literature of the past few years.
In this regard, the farthest-reaching insight 	
seems to have been achieved for systems in which signal-dependent migration mechanisms of the form addressed 
in (\ref{0}) are coupled to equations reflecting signal {\em production} through individuals, rather than consumption as in (\ref{0}). 
Indeed, for various classes of the key ingredient $\phi$, significant progress 
could be achieved for corresponding initial-boundary value problem associated with
\be{01}
	\left\{ \begin{array}{ll}
	u_t = \Delta \big(u\phi(v)\big),
	\qquad & x\in\Om, \ t>0, \\[1mm]
	\tau v_t = \Delta v-v+u,
	\qquad & x\in\Om, \ t>0, 
	\ear
\ee
both in the fully parabolic case when $\tau=1$, and in the simplified parabolic-elliptic version obtained on letting $\tau=0$.
Besides basic results on global solvability (see \cite{ahn_yoon}, \cite{jiang_laurencot}, \cite{fujie_senba}, 
\cite{fujie_jiang_ACAP2021}, \cite{desvillettes}, \cite{jin_kim_wang}, \cite{wang_wang}, 
\cite{taowin_M3AS} and also \cite{win_NON}), 
several studies also include findings on asymptotic behavior, mainly identifying diffusion-dominated constellations 
in which solutions to either (\ref{01}) or certain closely related variants 
stabilize toward homogeneous equilibria (\cite{jiang_arxiv}, \cite{taowin_M3AS}, \cite{fujie_jiang_JDE2020}, 
\cite{jin_kim_wang}, \cite{liu_xu}, \cite{wenbin_lv_EECT}, \cite{wenbin_lv_PROCA}, \cite{wenbin_lv_ZAMP}, \cite{yifu_wang}).
Beyond this, some interesting recent developments have provided rigorous evidence for a strong structure-supporting potential
of such models, as already predicted by some numerical experiments (\cite{desvillettes}),
by revealing the occurrence of infinite-time blow-up phenomena (\cite{fujie_jiang_arxiv}).\abs
{\bf Motility degeneracies at small signal densities.} \quad
According to manifest positivity properties of $v$ inherent to the production-determined signal evolution mechanism therein
(\cite{fujie_diss}),	
considerations related to the bahavior of $\phi$ close to the origin seem of secondary relevance in the context of (\ref{01});
in particular, degeneracies due to either singular or vanishing asymptotics of $\phi$ near $v=0$ appear to have no significant
effects on corresponding solution theories, and thus have partially even been explicitly included in precedent analysis
of (\ref{0}) (\cite{ahn_yoon}, \cite{jiang_laurencot}, \cite{fujie_jiang_ACAP2021}, \cite{win_NON}).\abs
In stark contrast to this, addressing the taxis-consumption problem
(\ref{0}) with its evident tendency toward enhancing small signal densities seems to require a distinct focus
on the behavior of $\phi(v)$ for small values of $v$, where especially migration-limiting mechanisms appear relevant
in the modeling of bacterial motion on nutrient-poor envoronments (\cite{kawasaki_JTB1997}, \cite{plaza}).
In line with this, the present manuscript will be concerned with (\ref{0}) under the assumption that $\phi$ be small near $v=0$,
while otherwise being fairly arbitrary for large $v$ and hence possibly retaining essential decay features at large signal
densities that have underlain the modeling hypotheses in \cite{fu} and \cite{liu}. \abs
Already at the level of questions related to mere solvability, this degenerate framework seems to bring about noticeable challenges,
especially in the application-relevant case when the initial signal distribution is small, or when $v_0$ even attains zeroes.
Adequately coping with such degeneracies in the course of an existence analysis for (\ref{0}) will accordingly form our first 
objective, and our attempt to accomplish this will lead us to a more general problem from basic parabolic theory.\abs
{\bf Quantifying positivity in a linear parabolic problem. Main results I.} \quad
Specifically, a crucial step in our approach will, guided by the ambition to make efficient use of the dissipative action in the 
first sub-problem of (\ref{0}), 
consist in establishing appropriate lower bounds
for the corresponding second components of solutions $(\ueps,\veps)$ to suitably regularized variants of (\ref{0}) (cf.~(\ref{0eps})).
Inter alia due to this approximation-based procedure, the necessity to thus simultaneously deal with whole solution {\em families},
instead of just one single object, seems to restrict accessibility to classical strong maximum principles.\abs
A crucial part of our analysis will accordingly be concerned with the derivation of positive pointwise lower 
bounds for families of nonnegative solutions $V$ to $V_t=\Del V - U(x,t)V$ possibly attaining zeroes initially,
under adequate assumptions on $U$ which are mild enough so as to allow for an application to $V:=\veps$ and $U:=\ueps$
on the basis of a priori information on $\ueps$ that can separately be obtained (cf.~the further discussion near (\ref{fi}) below).\abs
To address this in a context conveniently general, as an object of potentially independent interest we shall examine 
this question for the problem
\be{6.3}
	\lball
	V_t = \Del V + \na \cdot \big(a(x,t)V\big) + b(x,t) V,
	\qquad & x\in \Om, \ t\in (0,T_0), \\[1mm]
	\frac{\pa V}{\pa\nu}=0,
	\qquad & x\in\pO, \ t\in (0,T_0), \\[1mm]
	V(x,0)=V_0(x),
	\qquad & x\in\Om,
	\ear
\ee
and our main result in this direction will indeed establish a quantitative link between basic properties of $V_0$ as well as
$a$ and $b$ on the one hand, and positivity features of $V$ on the other:
\newcommand{\pp}{\Pi}
\begin{prop}\label{prop6}
  Let $n\ge 1$ and $\Om\subset\R^n$ be a bounded domain with smooth boundary,
  and suppose that $p_1\ge 2$, $p_2\ge 1$, 
  $q_1>2$ and $q_2>1$
  are such that
  \be{6.01}
	\frac{1}{q_1} + \frac{n}{2p_1} < \frac{1}{2}
	\qquad \mbox{and} \qquad
	\frac{1}{q_2} + \frac{n}{2p_2} < 1.
  \ee
  Then given any $L>0$, $T>0$ and $\tau\in (0,T)$ one can find $C(p_1,p_2,q_1,q_2,L,T,\tau)>0$
  with the property that whenever $T_0\in (0,T]$ as well as $a\in C^{1,0}(\bom\times [0,T_0);\R^n)$, 
  $b\in C^0(\bom\times [0,T_0))$ and
  $V\in C^0(\bom\times [0,T_0)) \cap C^{2,1}(\bom\times (0,T_0))$ are such that $a\cdot\nu=0$ on $\pO\times (0,T_0)$, that
  \be{6.1}
	\int_0^{T_0} \|a(\cdot,t)\|_{L^{p_1}(\Om)}^{q_1} \le L
	\qquad \mbox{and} \qquad
	\int_0^{T_0} \|b(\cdot,t)\|_{L^{p_2}(\Om)}^{q_2} \le L,
  \ee
  that
  \be{6.2}
	0 \le V_0 \le L
	\quad \mbox{in $\Omega$}
	\qquad
	\mbox{and} \qquad
	\io V_0 \ge \frac{1}{L},
  \ee
  and that (\ref{6.3}) holds,
  we have
  \be{6.4}
	V(x,t) \ge C(p_1,p_2,q_1,q_2,L,T,\tau)
	\qquad \mbox{for all $x\in\Om$ and } t\in (\tau,T_0).
  \ee
\end{prop}
Here we note that the hypotheses in (\ref{6.1}) and (\ref{6.01}) cannot be substantially relaxed, not even in the simple
case when $a\equiv 0$:
\begin{prop}\label{prop99}
  Let $n\ge 1$ and $\Om\subset\R^n$ be a bounded domain with smooth boundary,
  let $p\ge 1$ and $q\ge 1$ be such that
  \be{99.1}
	\frac{1}{q} + \frac{n}{2p}>1,
  \ee
  and let $T>0$ and $x_0\in\Om$.
  Then there exist $L>0$, $(b_k)_{k\in\N} \subset C^\infty(\bom\times [0,T])$ and a positive function
  $V_0\in C^\infty(\bom)$ such that
  \be{99.2}
	\int_0^T \|b_k(\cdot,t)\|_{L^p(\Om)}^q dt \le L
	\qquad \mbox{for all } k\in\N
  \ee
  and that (\ref{6.2}) holds, but that for each $k\in\N$ one can find $V_k\in C^0(\bom\times [0,T]) \cap C^{2,1}(\bom\times (0,T))$
  such that
  \be{99.3}
	\lball
	V_{kt} = \Del V_k + b_k(x,t) V_k,
	\qquad & x\in \Om, \ t\in (0,T), \\[1mm]
	\frac{\pa V_k}{\pa\nu}=0,
	\qquad & x\in\pO, \ t\in (0,T), \\[1mm]
	V_k(x,0)=V_0(x),
	\qquad & x\in\Om,
	\ear
  \ee
  and that
  \be{99.4}
	V_k(x_0,t) \to 0
	\qquad \mbox{as } k\to\infty.
  \ee
\end{prop}
{\bf Global solvability and large time behavior in (\ref{0}). Main results II.} \quad
In line with the requirements expressed in (\ref{6.1}) and (\ref{6.01}), our application of Proposition \ref{prop6}
to approximate solutions of (\ref{0}) needs to be preceded by the derivation of suitable integral bounds for
the respective first components $\ueps$.
When addressing this in the setting
of a standard $L^p$ testing procedure, we will be forced to appropriately control 	
ill-signed cross-diffusive contributions by means of correspondingly signal-weighted and hence weakened dissipation rates
(cf.~(\ref{5.22})).
This will be achieved in the course of a further essential step in our analysis, to be accomplished in Lemma \ref{lem5}, by utilizing 
a functional inequality of the form
\bea{fi}
	\io \frac{\vp^p}{\psi} |\na\psi|^2
	&\le& \eta \io \vp^{p-2} \psi |\na\vp|^2
	+ \eta \io \vp\psi \nn\\
	& & + C(p) \cdot \Big(1+\frac{1}{\eta}\Big) \cdot \Bigg\{ \io \vp^p + \bigg\{ \io \vp \bigg\}^{2p-1} \Bigg\}
		\cdot \io \frac{|\na\psi|^4}{\psi^3},
\eea
to be deduced in Lemma \ref{lem4}
for smooth $\vp\ge 0$ and $\psi>0$ and any $p\ge 2$ and $\eta>0$ with some $C(p)>0$ in one- and two-dimensional domains.
This will enable us to establish $L^p$ bounds for $\ueps$ in actually any $L^p$ space with finite $p$, 
and a subsequent application of Proposition \ref{prop6}, as thereby facilitated, will thereupon provide accessibility
to arguments from well-established parabolic regularity theories so as to finally yield $C^{2+\theta,1+\frac{\theta}{2}}$ estimates
within the range where said positivity result holds, that is, locally away from the temporal origin
(Lemma \ref{lem9} and Lemma \ref{lem11}).\abs
In consequence, this will enable us to establish the following result on global solvability of (\ref{0}) by functions
which are even smooth for all positive times, provided that $\phi$ and the initial data comply with mild assumptions which inter alia
allow for large classes of merely nonnegative $v_0$:
\begin{theo}\label{theo17}
  Let $n\in \{1,2\}$ and $\Om\subset\R^n$ be a bounded convex domain with smooth boundary, assume that
  \be{phi}
	\phi\in C^1([0,\infty)) \cap C^3((0,\infty))
	\quad
	\mbox{is such that $\phi(0)=0, \phi'(0)>0$ and $\phi>0$ on $(0,\infty)$,}
  \ee
  and suppose that
  \be{init}
	\lbal
	u_0\in W^{1,\infty}(\Om)
	\mbox{ is nonnegative with $u_0\not\equiv 0$, \quad and that} \\[1mm]
	v_0\in W^{1,\infty}(\Om)
	\mbox{ is nonnegative with $v_0\not\equiv 0$ and $\sqrt{v_0} \in W^{1,2}(\Om)$.}
	\ear
  \ee
  Then there exist functions
  \be{17.1}
	\lbal
	u\in C^{2,1}(\bom\times (0,\infty)) 
	\qquad \mbox{and} \\[1mm]
	v\in C^0(\bom\times [0,\infty)) \cap C^{2,1}(\bom\times (0,\infty)) 
	\ear
  \ee
  such that $u>0$ and $v>0$ in $\bom\times (0,\infty)$, and that $(u,v)$ solves (\ref{0}) in that in the classical pointwise
  sense we have $u_t=\Del \big(u\phi(v)\big)$ and $v_t=\Del v-uv$ in $\Om\times (0,\infty)$ 
  and $\frac{\pa u}{\pa\nu}=\frac{\pa v}{\pa\nu}=0$ on $\pO\times (0,\infty)$ as well as $v(\cdot,0)=v_0$ in $\Om$,
  and that
  \be{17.2}
	u(\cdot,t) \wto u_0
	\quad \mbox{in $L^p(\Om)$ for all $p\ge 1$ \quad as } t\searrow 0.
  \ee
  Moreover, this solution has the property that for each $p\ge 1$ there exists $C(p)>0$ fulfilling
  \bas
	\|u(\cdot,t)\|_{L^p(\Om)} + \|v(\cdot,t)\|_{W^{1,\infty}(\Om)} \le C(p)
	\qquad \mbox{for all } t>0.
  \eas
\end{theo}
Next focusing on the qualitative behavior of the solutions gained above, we shall 
make use of the decay information contained in an inequality of the form
\be{el}
	\int_0^\infty \io uv \le \io v_0,
\ee
as constituting one of the most elementary features of the second equation in (\ref{0}), to assert a bound in the style
of an estimate in $BV([0,\infty);(W_N^{2,\infty}(\Om))^\star)$ for $u$, where
$W_N^{2,\infty}(\Om):=\big\{\vp\in W^{2,\infty}(\Om) \ | \ \frac{\pa\vp}{\pa\nu}=0 \mbox{ on } \pO \big\}$
(Lemma \ref{lem14}).
Through interpolation, this will imply the essential part of the following result on large time stabilization of each among the
solutions obtained in Theorem \ref{theo17}:
\begin{theo}\label{theo18}
  Let $n\in \{1,2\}$ and $\Om\subset\R^n$ be a bounded convex domain with smooth boundary, and assume (\ref{phi}) 
  as well as (\ref{init}). Then there exists a nonnegative function $u_\infty\in \bigcap_{p\ge 1} L^p(\Om)$
  such that $\io u_\infty=\io u_0$, and that as $t\to\infty$, the solution $(u,v)$ of (\ref{0}) from Theorem \ref{theo17} satisfies
  \be{18.1}
	u(\cdot,t) \wto u_\infty
	\qquad \mbox{in $L^p(\Om)$ for all $p\ge 1$}
  \ee
  and
  \be{18.2}
	v(\cdot,t) \wsto 0
	\qquad \mbox{in } W^{1,\infty}(\Om).
  \ee
\end{theo}
As a natural question related to the latter result, we finally address the problem of describing the limit functions
$u_\infty$ appearing in (\ref{18.1}). To put this in perspective, let us recall that the literature has identified numerous
situations in which when accompanied by {\em non-degenerate} diffusion, 
taxis-type cross-diffusive interaction with absorptive signal evolution mechanisms as in (\ref{0}) 
leads to asymptotic prevalence of spatial homogeneity:
Indeed, not only (\ref{0}) with strictly positive $\phi$ (\cite{liwin2}), but also
a considerable variety of chemotaxis-consumption systems has been shown to have the common feature that
for widely arbitrary initial data, corresponding solutions stabilize toward {\em constant} states in their first component
(cf.~\cite{taowin_consumption}, \cite{lankeit_M3AS}, \cite{win_JDE_homogeneity} for some small selection of examples).\abs
A noticeable difference to this type of behavior, and hence a qualitative effect genuinely due to the diffusion
degeneracy in (\ref{0}), will be revealed in our final result: By making appropriate use of the quantitative 
information contained in (\ref{el}), we can derive a criterion, in its essence 
apparently reflecting quite well the nutrient-poor situation relevant to applications (\cite{kawasaki_JTB1997}, \cite{plaza}),
for the limit function in (\ref{18.1}) to be {\em nonconstant}:
\begin{theo}\label{theo19}
  Let $n\in \{1,2\}$ and $\Om\subset\R^n$ be a bounded convex domain with smooth boundary, suppose that (\ref{phi}) holds, and
  let $u_0\in W^{1,\infty}(\Om)$ be nonnegative with $u_0\not\equiv const.$
  Then for all $K>0$ there exists $\delta(K)>0$ with the property that whenever $v_0\in W^{1,\infty}(\Om)$ is nonnegative
  with $\sqrt{v_0}\in W^{1,2}(\Om)$ and such that
  \bas
	0<\|v_0\|_{L^\infty(\Om)}\le K
	\qquad \mbox{and} \qquad
	\io v_0 \le \delta(K),
  \eas
  the corresponding limit function obtained in Theorem \ref{theo18} satisfies
  \bas
	u_\infty \not\equiv const.
  \eas
\end{theo}
\mysection{A quantitative strong maximum principle. Proof of Proposition \ref{prop6}}
Let us first turn our attention to the most essential among our tools, by namely focusing on the positivity property claimed
in Proposition \ref{prop6}. 
Our reasoning in this direction will at its core be based on a comparison argument applied to the function $W:=\ln \frac{C}{V}$
which for suitably large $C$ depending on the parameters in Proposition \ref{prop6}, given any $V$ fulfilling (\ref{6.3}) indeed
satisfies an inhomogeneous linear parabolic inequality (cf.~(\ref{6.76})).
As an essential preparation for this, our derivation of some quantitative information on immediate smoothing of $W$ into
$L^1(\Om)$ will rely on a Poincar\'e type inequality, applicable here thanks to a short-time positive lower bound for $\io V$
due to (\ref{6.1}) (see (\ref{6.10})), which facilitates to make appropriate use of a first-order superlinear absorptive 
contribution to the evolution of $\io \ln \frac{\delta}{V}$ for suitably chosen $\delta>0$ (see (\ref{6.777})).\abs
Through this type of design, our strategy is able to cope with the mild regularity requirements in Proposition \ref{prop6},
and thereby, unlike alternative approaches based on lower estimates for Green's functions (\cite{choulli_kayser}),
especially remains applicable throughout the essentially optimal parameter range described by (\ref{6.01});
in particular, for our subsequently performed analysis of (\ref{0}) it will be of crucial importance that 
our argument in Proposition \ref{prop6} is robust enough so as to make do without requiring $L^\infty$ bounds for $b$.\abs
\proofc of Proposition \ref{prop6}. \quad
  We abbreviate $\pp:=(p_1,q_1,p_2,q_2)$ and 
  first recall known regularization features of the Neumann heat semigroup $(e^{t\Del})_{t\ge 0}$ on $\Om$ 
  (\cite{win_JDE2010}, \cite{FIWY}) to fix positive constants $c_1(\pp,T)$, $c_2(\pp,T)$, $c_3(\pp,T)$ and $c_4(T)>0$ 
  such that for any $t\in (0,T)$,
  \be{6.5}
	\|e^{t\Del} \na\cdot \vp\|_{L^\infty(\Om)} \le c_1(\pp,T) t^{-\frac{1}{2}-\frac{n}{2p_1}} \|\vp\|_{L^{p_1}(\Om)}
	\qquad \mbox{for all $\vp\in C^1(\bom;\R^n)$ fulfilling $\vp\cdot\nu|_{\pO}=0$}
  \ee
  and
  \be{6.56}
	\|e^{t\Del} \vp\|_{L^\infty(\Om)} \le c_2(\pp,T) t^{-\frac{n}{p_1}} \|\vp\|_{L^\frac{p_1}{2}(\Om)}
	\qquad \mbox{for all } \vp\in C^0(\bom)
  \ee
  and
  \be{6.55}
	\|e^{t\Del} \vp\|_{L^\infty(\Om)} \le c_3(\pp,T) t^{-\frac{n}{2p_2}} \|\vp\|_{L^{p_2}(\Om)}
	\qquad \mbox{for all } \vp\in C^0(\bom)
  \ee
  as well as
  \be{6.6}
	\|e^{t\Del} \vp\|_{L^\infty(\Om)} \le c_4(T) t^{-\frac{n}{2}} \|\vp\|_{L^1(\Om)}
	\qquad \mbox{for all } \vp\in C^0(\bom).
  \ee
  Using that 
  \be{6.61}
	\big(\frac{1}{2}+\frac{n}{2p_1}\Big)\cdot\frac{q_1}{q_1-1} <1
	\qquad \mbox{and} \qquad
	\frac{n}{2p_2} \cdot \frac{q_2}{q_2-1}<1
  \ee
  according to (\ref{6.01}), we can thereafter rely on Beppo Levi's theorem
  to fix $\mu(\pp,L,T)>0$ suitably large such that
  \be{6.66}
	c_1(\pp,T) L^\frac{1}{q_1} \cdot \bigg\{ 
		\int_0^T \sigma^{-(\frac{1}{2}+\frac{n}{2p_1})\cdot\frac{q_1}{q_1-1}} 
		e^{-\frac{\mu(\pp,L,T) q_1}{q_1-1} \cdot \sigma} d\sigma
		\bigg\}^\frac{q_1-1}{q_1} \le \frac{1}{4}
  \ee
  and
  \be{6.67}
	c_3(\pp,T) L^\frac{1}{q_2} \cdot \bigg\{ 
		\int_0^T \sigma^{-\frac{n}{2p_2} \cdot \frac{q_2}{q_2-1}} e^{-\frac{\mu(\pp,L,T) q_2}{q_2-1}\cdot \sigma} d\sigma
		\bigg\}^\frac{q_2-1}{q_2} \le \frac{1}{4}.
  \ee
  We moreover employ a consequence of a Poincar\'e type inequality (\cite[Lemma 8.4 and appendix]{lankeit_win_NoDEA},
  \cite[Lemma 4.3]{taowin_persistence}) to choose $c_5(\pp,L,T)>0$ in such a way that whenever $\delta>0$,
  \bea{6.7}
	& & \hs{-20mm}
	\frac{1}{2} \io \frac{|\na\vp|^2}{\vp^2}
	\ge c_5(\pp,L,T) \cdot \bigg\{ \io \ln \frac{\delta}{\vp} \bigg\}_+^2 \nn\\
	& & \hs{0mm}
	\qquad \mbox{for all $\vp\in C^1(\bom)$ such that $\vp>0$ in $\bom$ and \quad}
	\big| \{\vp>\delta\} \big| \ge \frac{1}{4Lc_6(\pp,L,T)},
  \eea
  where
  \be{6.77}
	c_6(\pp,L,T):=2L e^{\mu(\pp,L,T) T}.
  \ee
  We now suppose that $T_0\in (0,T]$ and that $a, V_0$ and $V$ have the listed properties, and begin our derivation of (\ref{6.4})
  by relying on a Duhamel representation associated with (\ref{6.3}) to see that thanks to the maximum principle, the identity
  $a\cdot\nu|_{\pO\times (0,T_0)}=0$,
  (\ref{6.5}), (\ref{6.55}) and the H\"older inequality, the continuous function $y$ given by 
  $y(t):=e^{-\mu(\pp,L,T) t} \|V(\cdot,t)\|_{L^\infty(\Om)}$, $t\in [0,T_0)$, satisfies
{
\allowdisplaybreaks
  \bas
	y(t)
	&=& e^{-\mu(\pp,L,T) t} 
	\bigg\| e^{t\Del} V_0 + \int_0^t e^{(t-s)\Del} \na\cdot \big\{ a(\cdot,s)V(\cdot,s)\big\} ds 
	+ \int_0^t e^{(t-s)\Del} \big\{ b(\cdot,s) V(\cdot,s) \big\} ds \bigg\|_{L^\infty(\Om)} \\
	&\le& e^{-\mu(\pp,L,T) t} \|V_0\|_{L^\infty(\Om)}
	+ c_1(\pp,T) e^{-\mu(\pp,L,T) t} 
	\int_0^t (t-s)^{-\frac{1}{2}-\frac{n}{2p_1}} \big\| a(\cdot,s) V(\cdot,s)\big\|_{L^{p_1}(\Om)} ds \\
	& & + c_3(\pp,T) e^{-\mu(\pp,L,T) t}
	\int_0^t (t-s)^{-\frac{n}{2p_2}} \big\| b(\cdot,s) V(\cdot,s)\big\|_{L^{p_2}(\Om)} ds \\
	&\le& e^{-\mu(\pp,L,T) t} \|V_0\|_{L^\infty(\Om)}
	+ c_1(\pp,T) e^{-\mu(\pp,L,T) t} 
		\int_0^t (t-s)^{-\frac{1}{2}-\frac{n}{2p_1}} \|a(\cdot,s)\|_{L^{p_1}(\Om)} \|V(\cdot,s)\|_{L^\infty(\Om)} ds \\
	& & + c_3(\pp,T) e^{-\mu(\pp,L,T) t}
		\int_0^t (t-s)^{-\frac{n}{2p_2}} \|b(\cdot,s)\|_{L^{p_2}(\Om)} \|V(\cdot,s)\|_{L^\infty(\Om)} ds \\
	&\le& e^{-\mu(\pp,L,T) t} \|V_0\|_{L^\infty(\Om)} \\
	& & + c_1(\pp,T) \cdot \bigg\{ \int_0^t \|a(\cdot,s)\|_{L^{p_1}(\Om)}^{q_1} ds \bigg\}^\frac{1}{q_1} \times \\
	& & \hs{20mm}
		\times \bigg\{ \int_0^t (t-s)^{-(\frac{1}{2}+\frac{n}{2p_1}) \cdot \frac{q_1}{q_1-1}} 
		e^{- \frac{\mu(\pp,L,T) q_1}{q_1-1} \cdot (t-s)} ds \bigg\}^\frac{q_1-1}{q_1} 
		\cdot \|y\|_{L^\infty((0,t))} \\
	& & + c_3(\pp,T) \cdot \bigg\{ \int_0^t \|b(\cdot,s)\|_{L^{p_2}(\Om)}^{q_2} ds \bigg\}^\frac{1}{q_2} \times \\
	& & \hs{20mm}
		\times \bigg\{ \int_0^t (t-s)^{-\frac{n}{2p_2} \cdot \frac{q_2}{q_2-1}} 
		e^{- \frac{\mu(\pp,L,T) q_2}{q_2-1} \cdot (t-s)} ds \bigg\}^\frac{q_2-1}{q_2} 
		\cdot \|y\|_{L^\infty((0,t))} 
  \eas
}
  for all $t\in (0,T_0)$.
  Therefore, (\ref{6.2}) and (\ref{6.1}) together with (\ref{6.66}) and (\ref{6.67}) ensure that
  \bas
	y(t) 
	&\le& L 
	+ c_1(\pp,T) L^\frac{1}{q_1} \cdot \bigg\{ 
		\int_0^{T_0} \sigma^{-(\frac{1}{2}+\frac{n}{2p_1})\cdot\frac{q_1}{q_1-1}} 
		e^{-\frac{\mu(\pp,L,T) q_1}{q_1-1} \cdot \sigma} d\sigma
		\bigg\}^\frac{q_1-1}{q_1} 
	\cdot \|y\|_{L^\infty((0,t))} \\
	& & + c_3(\pp,T) L^\frac{1}{q_2} \cdot \bigg\{ 
		\int_0^{T_0} \sigma^{-\frac{n}{2p_2} \cdot \frac{q_2}{q_2-1}} e^{-\frac{\mu(\pp,L,T) q_2}{q_2-1} \cdot \sigma} d\sigma
		\bigg\}^\frac{q_2-1}{q_2} 
	\cdot \|y\|_{L^\infty((0,t))} \\
	&\le& L + \frac{1}{4} \|y\|_{L^\infty((0,t))} + \frac{1}{4} \|y\|_{L^\infty((0,t))}
	\qquad \mbox{for all } t\in (0,T_0),
  \eas
  from which it follows that, in line with (\ref{6.77}),
  \be{6.666}
	\|V(\cdot,t)\|_{L^\infty(\Om)} \le c_6(\pp,L,T)
	\qquad \mbox{for all } t\in (0,T_0).
  \ee
  Again since $a\cdot\nu=0$ on $\pO\times (0,T_0)$, in view of the H\"older inequality this especially ensures that
  \bas
	\frac{d}{dt} \io V
	&=& \io b(x,t) V \\
	&\ge& - c_6(\pp,L,T) |\Om|^\frac{p_2-1}{p_2} \|b(\cdot,t)\|_{L^{p_2}(\Om)} 
	\qquad \mbox{for all } t\in (0,T_0)
  \eas
  and that hence, by (\ref{6.2}) and (\ref{6.1}),
  \bea{6.10}
	\io V(\cdot,t)
	&\ge& \io V_0 - c_6(\pp,L,T) |\Om|^\frac{p_2-1}{p_2} \int_0^t \|b(\cdot,s)\|_{L^{p_2}(\Om)} ds \nn\\
	&\ge& \frac{1}{L} - c_6(\pp,L,T) |\Om|^\frac{p_2-1}{p_2} \cdot 
		\bigg\{ \int_0^t \|b(\cdot,s)\|_{L^{p_2}(\Om)}^{q_2} ds \bigg\}^\frac{1}{q_2} \cdot t^\frac{q_2-1}{q_2} \nn\\
	&\ge& \frac{1}{L} - c_6(\pp,L,T) |\Om|^\frac{p_2-1}{p_2} L^\frac{1}{q_2} t^\frac{q_2-1}{q_2} \nn\\
	&\ge& \frac{1}{2L}
	\qquad \mbox{for all } t\in (0,\wh{t}_1),
  \eea  
  where 
  \be{6.100}
	\wh{t}_1:=\min \big\{ t_1,T_0 \big\}
	\qquad \mbox{with} \qquad
	t_1\equiv t_1(\pp,L,T):= \Big\{ 2c_6(\pp,L,T) |\Om|^\frac{p_2-1}{p_2} L^\frac{q_2+1}{q_2} \Big\}^{-\frac{q_2}{q_2-1}}.
  \ee
  Combining this with (\ref{6.666}), we see that for
  $\delta(L):=\frac{1}{4|\Om| L}$ we have
  \bas
	\frac{1}{2L}
	&\le& \int_{\{V(\cdot,t)\le\delta(L)\}} V(\cdot,t)
	+ \int_{\{V(\cdot,t)>\delta(L)\}} V(\cdot,t) \\
	&\le& \delta(L) |\Om|
	+ c_6(\pp,L,T) \cdot \big| \{V(\cdot,t)>\delta(L)\} \big| \\
	&=& \frac{1}{4L} 
	+ c_6(\pp,L,T) \cdot \big| \{V(\cdot,t)>\delta(L)\} \big|
	\qquad \mbox{for all } t\in (0,\wh{t}_1)
  \eas
  and thus
  \bas
	\big|\{V(\cdot,t)>\delta(L)\} \big| \ge \frac{1}{4Lc_6(\pp,L,T)}
	\qquad \mbox{for all } t\in (0,\wh{t}_1).
  \eas
  We may therefore draw on (\ref{6.7}) to find that in the identity
  \be{6.777}
	\frac{d}{dt} \io \ln \frac{\delta(L)}{V}
	= - \io \frac{V_t}{V}
	= - \io \frac{|\na V|^2}{V^2}
	- \io a \cdot \frac{\na V}{V}
	- \io b,
  \ee
  valid throughout $(0,T_0)$ since clearly $V$ is positive on $\bom\times (0,T_0)$ by (\ref{6.2}) and the classical
  strong maximum principle, and again since $a\cdot\nu=0$ on $\pO\times (0,T_0)$, we can estimate
  \bas
	\frac{1}{2} \io \frac{|\na V|^2}{V^2}
	\ge c_5(\pp,L,T) \cdot \bigg\{ \io \ln \frac{\delta(L)}{V} \bigg\}_+^2
	\qquad \mbox{for all } t\in (0,\wh{t}_1).
  \eas
  As moreover 
  \bas
	- \io a\cdot\frac{\na V}{V} 
	&\le& \frac{1}{2} \io \frac{|\na V|^2}{V^2} 
	+ \frac{1}{2} \io |a|^2 \\
	&\le& \frac{1}{2} \io \frac{|\na V|^2}{V^2} 
	+ \frac{1}{2} |\Om|^\frac{p_1-2}{p_1} \|a(\cdot,t)\|_{L^{p_1}(\Om)}^2 
	\qquad \mbox{for all } t\in (0,T_0)
  \eas
  and
  \bas
	- \io b \le |\Om|^\frac{p_2-1}{p_2} \|b(\cdot,t)\|_{L^{p_2}(\Om)}
	\qquad \mbox{for all } t\in (0,T_0)
  \eas
  by the H\"older inequality, this implies that if we let 
  $c_7(\pp):=\max \big\{ \frac{1}{2} |\Om|^\frac{p_1-2}{p_1} \, , \, |\Om|^\frac{p_2-1}{p_2} \big\}$, then
  $z(t):=\io \ln \frac{\delta(L)}{V(\cdot,t)}$,
  $t\in (0,\wh{t}_1)$, has the property that
  \bas
	z'(t) \le -c_5(\pp,L,T) z_+^2(t) + c_7(\pp) \|a(\cdot,t)\|_{L^{p_1}(\Om)}^2 + c_7(\pp) \|b(\cdot,t)\|_{L^{p_2}(\Om)}
	\qquad \mbox{for all } t\in (0,\wh{t}_1).
  \eas
  By means of an ODE comparison argument, this can be seen to entail that with
  \be{6.11}
	h(t):=c_7(\pp) \int_0^t \|a(\cdot,s)\|_{L^{p_1}(\Om)}^2 ds
	+ c_7(\pp) \int_0^t \|b(\cdot,s)\|_{L^{p_2}(\Om)} ds,
	\qquad t\in (0,T_0),
  \ee
  we have
  \be{6.12}
	z(t) \le \frac{1}{c_5(\pp,L,T) t} + h(t)
	\qquad \mbox{for all } t\in (0,\wh{t}_1),
  \ee
  because for each $\eta\in (0,\wh{t}_1)$,
  \bas
	\ov{z}(t):=\frac{1}{c_5(\pp,L,T)\cdot (t-\eta)} + h(t),
	\qquad t>\eta,
  \eas
  satisfies
  \bas
	& & \hs{-20mm}
	\ov{z}'(t) + c_5(\pp,L,T) \ov{z}_+^2(t) 
	- c_7(\pp) \|a(\cdot,t)\|_{L^{p_1}(\Om)}^2 - c_7(\pp) \|b(\cdot,t)\|_{L^{p_2}(\Om)} \\
	&=& \bigg\{ -\frac{1}{c_5(\pp,L,T) \cdot (t-\eta)^2} + h'(t) \bigg\} 
	+ c_5(\pp,L,T) \cdot \bigg\{ \frac{1}{c_5(\pp,L,T)\cdot (t-\eta)} + h(t) \bigg\}^2 \\
	& & - c_7(\pp) \|a(\cdot,t)\|_{L^{p_1}(\Om)}^2 - c_7(\pp) \|b(\cdot,t)\|_{L^{p_2}(\Om)} \\
	&=& \frac{2h(t)}{t-\eta} + c_5(\pp,L,T) h^2(t) \\
	&\ge& 0
	\qquad \mbox{for all } t\in (\eta,\wh{t}_1)
  \eas
  according to (\ref{6.11}).
  In order to 
  make this applicable so as to accomplish the final step of our argument, we note that once more due to the H\"older inequality,
  (\ref{6.1}) ensures that
  \bas
	h(t) &\le& c_7(\pp) t^\frac{q_1-2}{q_1} \cdot \bigg\{ \int_0^t \|a(\cdot,s)\|_{L^{p_1}(\Om)}^{q_1} ds \bigg\}^\frac{2}{q_1}
	+ c_7(\pp) t^\frac{q_2-1}{q_2} \cdot \bigg\{ \int_0^t \|b(\cdot,s)\|_{L^{p_2}(\Om)}^{q_2} ds \bigg\}^\frac{1}{q_2} \\
	&\le& c_8(\pp,L,T):= c_7(\pp) T^\frac{q_1-2}{q_1} L^\frac{2}{q_1}
	+ c_7(\pp) T^\frac{q_2-1}{q_2} L^\frac{1}{q_2}
	\qquad \mbox{for all } t\in (0,T_0),
  \eas
  and that thus the function $W$ defined by
  \bas
	W(x,t):=\ln \frac{c_6(\pp,L,T)}{V(x,t)},
	\qquad x\in\bom, \ t\in (0,T_0),
  \eas
  nonnegative throughout $\bom\times (0,T_0)$ thanks to (\ref{6.666}), satisfies
  \bea{6.13}
	\|W(\cdot,t)\|_{L^1(\Om)}
	&=& \io \ln \Big\{ \frac{\delta(L)}{V(\cdot,t)} \cdot \frac{c_6(\pp,L,T)}{\delta(L)} \Big\} \nn\\[2mm]
	&\le& z(t) + \frac{|\Om| c_6(\pp,L,T)}{\delta(L)} \nn\\[2mm]
	&\le& \frac{1}{c_5(\pp,L,T) t} + c_9(\pp,L,T)
	\qquad \mbox{for all } t\in (0,\wh{t}_1)
  \eea
  with $c_9(\pp,L,T):=c_8(\pp,L,T) + \frac{|\Om| c_6(\pp,L,T)}{\delta(L)}$.
  To derive (\ref{6.4})
  from this, we let $\tau\in (0,T)$ be given and note that we only need to consider the case when $T_0>\tau$,
  in which (\ref{6.100}) warrants that
  $\wh{t}_1 \ge t_2 \equiv t_2(\pp,L,T,\tau):=\min\{ t_1,\tau\}$.
  As (\ref{6.3}) together with Young's inequality implies that
  \bea{6.76}
	W_t &=& \Del W - |\na W|^2 - \frac{1}{V} \na\cdot\big( a(x,t) V\big) - b(x,t) \nn\\
	&=& \Del W - |\na W|^2 - \na \cdot a(x,t) - a(x,t)\cdot\na W - b(x,t) \nn\\
	&\le& \Del W - \na \cdot a(x,t) + \frac{1}{4} |a(x,t)|^2 - b(x,t)
	\qquad \mbox{in } \Om\times (0,T_0),
  \eea
  according to the comparison principle we may use (\ref{6.5}) and (\ref{6.55}) 
  now together with (\ref{6.56}) and (\ref{6.6}) to infer on the basis of a corresponding
  variation-of-constants representation that thanks to (\ref{6.13}), the H\"older inequality and (\ref{6.1}),
  \bas
	W(\cdot,t)
	&\le& e^{(t-\frac{t_2}{2})\Del} W\Big(\cdot,\frac{t_2}{2}\Big) \\
	& & - \int_\frac{t_2}{2}^t e^{(t-s)\Del} \na\cdot a(\cdot,s) ds 
	+ \frac{1}{4} \int_\frac{t_2}{2}^t e^{(t-s)\Del} |a(\cdot,s)|^2 ds
	- \int_\frac{t_2}{2}^t e^{(t-s)\Del} b(\cdot,s) ds \\
	&\le& c_4(T) \cdot \Big(t-\frac{t_2}{2}\Big)^{-\frac{n}{2}} \Big\| W\Big(\cdot,\frac{t_2}{2}\Big)\Big\|_{L^1(\Om)} \\
	& & + c_1(\pp,T) \int_\frac{t_2}{2}^t (t-s)^{-\frac{1}{2}-\frac{n}{2p_1}} \|a(\cdot,s)\|_{L^{p_1}(\Om)} ds \\
	& & + \frac{c_2(\pp,T)}{4} \int_\frac{t_2}{2}^t (t-s)^{-\frac{n}{p_1}} \big\| |a(\cdot,s)|^2 \big\|_{L^\frac{p_1}{2}(\Om)} ds \\
	& & + c_3(\pp,T) \int_\frac{t_2}{2}^t (t-s)^{-\frac{n}{2p_2}} \|b(\cdot,s)\|_{L^{p_2}(\Om)} ds \\
	&\le& c_4(T) \cdot \Big(t-\frac{t_2}{2}\Big)^{-\frac{n}{2}} \cdot 
	\Big\{ \frac{2}{c_5(\pp,L,T) t_2} + c_9(\pp,L,T)\Big\}  \\[3mm]
	& & + c_{10}(\pp,L,T)
	\quad \mbox{in } \Om
	\qquad \mbox{for all } t\in \Big(\frac{t_2}{2},T_0\Big),
  \eas
  where 
  \bas
	c_{10}(\pp,L,T)
	&:=&
	c_1(\pp,L,T) L^\frac{1}{q_1} \cdot \bigg\{ \int_0^T \sigma^{-(\frac{1}{2}+ \frac{n}{2p_1})\cdot \frac{q_1}{q_1-1}} d\sigma
		\bigg\}^\frac{q_1-1}{q_1} \\
	& & + \frac{c_2(\pp,L,T)}{4} L^\frac{2}{q_1} \cdot \bigg\{ \int_0^T \sigma^{-\frac{n}{p_1} \cdot \frac{q_1}{q_1-2}} d\sigma
		\bigg\}^\frac{q_1-2}{q_1} \\
	& & + c_3(\pp,L,T) L^\frac{1}{q_2} \cdot \bigg\{ \int_0^T \sigma^{-\frac{n}{2p_2} \cdot \frac{q_2}{q_2-1}} d\sigma
		\bigg\}^\frac{q_2-1}{q_2}
  \eas
  is finite because of (\ref{6.61}), and of the fact that (\ref{6.01}) moreover warrants that 
  $\frac{n}{p_1} \cdot \frac{q_1}{q_1-2}<1$.
  Since $t_2\le\tau$, by definition of $W$ this particularly means that for all $x\in\Om$ and $t\in (\tau,T_0)$ we have
  \bas
	V(x,t)
	&\ge& c_6(p,L,T) \times \\
	& & \times \exp \bigg\{
	- c_4(T) \cdot \Big( \frac{t_2(\pp,L,T,\tau)}{2} \Big)^{-\frac{n}{2}} \cdot 
	\Big\{ \frac{2}{c_5(\pp,L,T) t_2(\pp,L,T,\tau)} + c_9(\pp,L,T)\Big\} \\
	& & \hs{15mm} 
	- c_{10}(\pp,L,T) \bigg\},
  \eas
  and that hence indeed (\ref{6.4}) holds with some $C(\pp,L,T,\tau)>0$ independent of $T_0$, $a, b, V_0$ and $V$.
\qed
Our construction of a counterexample in the case when instead of (\ref{6.01}) we have (\ref{99.1}) is much less involved:\abs
\proofc of Proposition \ref{prop99}. \quad
  We fix any $\alpha\in (0,1)$ and then use that $2n+\alpha-(1-\alpha)\xi^2 \to -\infty$ as $\xi\to\infty$ to pick a nonnegative 
  function $g\in C_0^\infty([0,\infty))$ such that
  \be{99.5}
	(\xi^2+1) g(\xi) \ge 2n+\alpha-(1-\alpha) \xi^2
	\qquad \mbox{for all } \xi\ge 0.
  \ee
  Without loss of generality assuming that $x_0=0$, we then choose $R>0$ and $R_0>R$ such that 
  $\ov{B}_R(0)\subset \Om \subset B_{R_0}(0)$, and for fixed $T>0$ taking $(T_k)_{k\in\N} \subset (T,T+1)$ such that $T_k\to T$
  as $k\to\infty$, we let
  \bas
	b_k(x,t):=-(T_k-t)^{-1} \cdot g\Big( (T_k-t)^{-\frac{1}{2}} |x|\Big),
	\qquad x\in\bom, \ t\in [0,T],
  \eas
  for $k\in\N$.
  Then since $T_k>T$, it follows that $b_k$ indeed belongs to $C^\infty(\bom\times [0,T])$ and, with $\omega_n:=n|B_1(0)|$, 
  due to the inclusion $\Om\subset B_{R_0}(0)$ satisfies
  \bas
	\io |b_k(x,t)|^p dx
	&\le& \omega_n \int_0^{R_0} r^{n-1} \cdot \Big\{ (T_k-t)^{-1} \cdot g\Big( (T_k-t)^{-\frac{1}{2}} r\Big) \Big\}^p dr \\
	&=& \omega_n \cdot (T_k-t)^{-p} \int_0^{R_0} r^{n-1} g^p\Big( (T_k-t)^{-\frac{1}{2}} r\Big) dr \\
	&=& \omega_n \cdot (T_k-t)^{\frac{n}{2}-p} \int_0^{(T_k-t)^{-\frac{1}{2}} R_0} \xi^{n-1} g^p(\xi) d\xi \\
	&\le& c_1 \cdot (T_k-t)^{\frac{n}{2}-p}
	\qquad \mbox{for all $t\in (0,T)$ and } k\in\N,
  \eas
  where $c_1:=\omega_n \int_0^\infty \xi^{n-1} g^p(\xi) d\xi$ is finite according to the boundedness of $\supp g$.
  Therefore,
  \bas
	\int_0^T \|b_k(\cdot,t)\|_{L^p(\Om)}^q dt
	\le c_1^\frac{q}{p} \int_0^T (T_k-t)^{(\frac{n}{2}-p)\cdot\frac{q}{p}} dt 
	\le c_1^\frac{q}{p} \int_0^{T_k} s^{(\frac{n}{2}-p)\cdot\frac{q}{p}} ds
	\qquad \mbox{for all } k\in\N,
  \eas
  so that since our assumption (\ref{99.1}) ensures that $(\frac{n}{2}-p)\cdot\frac{q}{p}>-1$, we can find $c_2>0$ fulfilling
  \be{99.6}
	\int_0^T \|b_k(\cdot,t)\|_{L^p(\Om)}^q dt
	\le c_2
	\qquad \mbox{for all } k\in\N;
  \ee
  writing $V_0(x):=\min\{ T^\alpha \, , \, (T+1)^{\alpha-1} R^2\}$, $x\in\bom$, we can thereupon fix $L>0$
  large enough such that besides (\ref{6.2}) we also have $L\ge c_2$.\abs
  It is then clear that thanks to the smoothness features of the constant function $V_0$ and of $(b_k)_{k\in\N}$, according
  to standard parabolic theory (\cite{LSU}) for any $k\in\N$ the problem (\ref{99.3}) admits a classical solution
  $V_k\in C^0(\bom\times [0,T]) \cap C^{2,1}(\bom\times (0,T))$ which, by nonpositivity of $b_k$ and the maximum principle, satisfies
  \be{99.7}
	V_k \le \min \big\{ T^\alpha \, , \, (T+1)^{\alpha-1} R^2 \big\}
	\qquad \mbox{in } \bom\times [0,T].
  \ee
  To derive (\ref{99.4}) from this, we let
  \bas
	\ov{V}_k(x,t):=(T_k-t)^\alpha f\Big( (T_k-t)^{-\frac{1}{2}} |x| \Big),
	\qquad (x,t)\in\bom\times [0,T], \ k\in\N,
  \eas
  with $f(\xi):=\xi^2+1$, $\xi\ge 0$,
  and use that $f'(\xi)=2\xi$ and $f''(\xi)=2$ for all $\xi\ge 0$ in verifying that for each $k\in\N$ and any 
  $(x,t)\in\Om\times (0,T)$, writing $\xi\equiv \xi(x,t;k):=(T_k-t)^{-\frac{1}{2}}|x|$ we have
  \bas
	\ov{V}_{kt} - \Del \ov{V}_k - b_k(x,t) \ov{V}_k
	&=& - \alpha (T_k-t)^{\alpha-1} f(\xi)
	+ \frac{1}{2} (T_k-t)^{\alpha-\frac{3}{2}} |x| f'(\xi) \\
	& & - (T_k-t)^{\alpha} \cdot \Big\{ (T_k-t)^{-1} f''(\xi) + \frac{n-1}{|x|} (T_k-t)^{-\frac{1}{2}} f'(\xi) \Big\}  \\
	& & + (T_k-t)^{-1} \cdot g(\xi) \cdot (T_k-t)^{\alpha} f(\xi) \\
	&=& (T_k-t)^{\alpha-1} \cdot \Big\{ -\alpha f(\xi) + \frac{\xi}{2} f'(\xi)
	- f''(\xi) - \frac{n-1}{\xi} f'(\xi) + g(\xi) f(\xi) \Big\} \\
	&=& (T_k-t)^{\alpha-1} \cdot \Big\{ (1-\alpha) \xi^2 - \alpha - 2n + g(\xi) \cdot (\xi^2+1) \Big\} \\[2mm]
	&\ge& 0
  \eas
  due to (\ref{99.5}).
  Since for any $x\in\pa B_R(0)$ and all $t\in (0,T)$ we have
  \bas
	\ov{V}_k(x,t)
	= (T_k-t)^{\alpha} \cdot \big\{ (T_k-t)^{-1} R^2 +1\big\}
	\ge (T_k-t)^{\alpha-1} R^2
	\ge (T+1)^{\alpha-1} R^2
	\ge V_k(x,t)
  \eas
  according to the inequalities $T_k<T+1$ and $\alpha<1$, and thanks to (\ref{99.7}), and since the latter moreover entails that
  \bas
	\ov{V}_k(x,0)
	\ge T_k^\alpha \cdot \big\{ T_k^{-1}  |x|^2 +1 \big\}
	\ge T_k^\alpha \ge T^\alpha \ge V_k(x,0)
	\qquad \mbox{for all } x\in B_R(0),
  \eas
  from the comparison principle we thus infer that $\ov{V}_k\ge V_k$ in $B_R(0) \times (0,T)$ for all $k\in\N$.
  Therefore, (\ref{99.4}) results upon observing that since $\alpha$ is positive,
  \bas
	\inf_{t\in (0,T)} \ov{V}_k(0,t)= (T_k-T)^\alpha \to 0
	\qquad \mbox{as } k\to\infty 
  \eas
  due to our requirement that $T_k\to T$ as $k\to\infty$.
\qed
\mysection{Analysis of (\ref{0}): $L^p$ bounds}
Our analysis of (\ref{0}) will now be launched by the observation that according to standard arguments from the theory
of Keller-Segel type cross-diffusion systems, for each $\eps\in (0,1)$ the regularized variant of (\ref{0}) given by
\be{0eps}
	\left\{ \begin{array}{ll}
	u_{\eps t} = \eps\Del\ueps + \Del \big(\ueps\phi(\veps)\big), 
	\qquad & x\in\Om, \ t>0, \\[1mm]
	v_{\eps t} =\Del \veps -\ueps\veps,
	\qquad & x\in\Om, \ t>0, \\[1mm]
	\frac{\partial \ueps}{\partial\nu}=\frac{\partial \veps}{\partial\nu}=0,
	\qquad & x\in\pO, \ t>0, \\[1mm]
	\ueps(x,0)=u_0(x), \quad \veps(x,0)=v_0(x),
	\qquad & x\in\Om.
	\end{array} \right.
\ee
admits local-in-time classical solutions enjoying a handy extensibility criterion:
\begin{lem}\label{lem_loc}
  Let $n\ge 1$ and $\Om\subset\R^n$ be a bounded domain with smooth boundary, and suppose that (\ref{phi}) and (\ref{init}) hold.
  Then for each $\eps\in (0,1)$ there exist $\tme\in (0,\infty]$ and functions
  \bas
	\lbal
	\ueps\in C^0(\bom\times [0,\tme)) \cap C^{2,1}(\bom\times (0,\tme)) \qquad \mbox{and} \\[1mm]
	\veps\in \bigcap_{q\ge 1} C^0([0,\tme); W^{1,q}(\Om)) \cap C^{2,1}(\bom\times (0,\tme)) 
	\ear
  \eas
  such that $\ueps\ge 0$ and $\veps>0$ in $\bom\times (0,\tme)$, that $(\ueps,\veps)$ solves (\ref{0eps}) in the classical sense
  in $\Om\times (0,\tme)$, and that
  \be{ext}
	\mbox{if $\tme<\infty$, \quad then \quad}
	\limsup_{t\nearrow\tme} \|\ueps(\cdot,t)\|_{L^\infty(\Om)}=\infty.
  \ee
  This solution satisfies
  \be{mass}
	\io \ueps(\cdot,t)=\io u_0
	\qquad \mbox{for all } t\in (0,\tme)
  \ee
  and
  \be{vinfty}
	\|\veps(\cdot,t)\|_{L^\infty(\Om)} \le \|v_0\|_{L^\infty(\Om)}
	\qquad \mbox{for all } t\in (0,\tme)
  \ee
  as well as
  \be{uv}
	\int_0^{\tme} \io \ueps\veps \le \io v_0.
  \ee
\end{lem}
\proof
  The statements on existence, positivity and extensibility 
  can be verified by following standard approaches in local existence theories of taxis-type parabolic systems
  (\cite{amann}, \cite{lankeit_singularconsumption}).		
  The mass conservation property immediately results from an integration in the first sub-problem in (\ref{0eps}),
  whereas (\ref{vinfty}) is a consequence of the comparison principle. Finally, the inequality (\ref{uv}) can be verified
  upon a time integration of the identity $\frac{d}{dt} \io \veps = - \io \ueps\veps$.
\qed
Throughout the sequel, unless otherwise stated we shall tacitly assume that (\ref{phi}) holds, and that
$n\le 2$ and $\Om\subset\R^n$ is a smoothly bounded convex domain, noting that the convexity requirement will be needed
from Lemma \ref{lem3} on, while the restriction on the spatial dimension will be relied on only in Lemma \ref{lem4} and its sequel.
Moreover, once $u_0$ and $v_0$ fulfilling (\ref{init}) have been fixed, by $(\ueps,\veps)$ and $\tme$ we shall exclusively
mean the objects provided by Lemma \ref{lem_loc}.\abs
For repeated later reference, let us explicitly state the following elementary implication of our assumptions on $\phi$,
and especially the requirement that $\phi'(0)$ be positive.
\begin{lem}\label{lem1}
  Let $K>0$. Then there exist $\lambda(K)>0$ and $\Lambda(K)>0$ such that if (\ref{init}) holds with $\|v_0\|_{L^\infty(\Om)} \le K$,
  we have
  \be{1.1}
	\lambda(K)\veps \le \phi(\veps) \le \Lambda(K)\veps
	\qquad \mbox{in } \Om\times (0,\tme)
  \ee
  and
  \be{1.2}
	|\phi'(\veps)| \le \Lambda(K)
	\qquad \mbox{in } \Om\times (0,\tme).
  \ee
\end{lem}
\proof
  Since $\phi(0)=0$, letting $\Lambda(K):=\|\phi'\|_{L^\infty((0,K))}$ we obtain that besides (\ref{1.2}), also the right
  inequality in (\ref{1.1}) holds due to (\ref{vinfty}). 
  For the same reason, the l'hospital rule ensures that $\rho(\xi):=\frac{\phi(\xi)}{\xi}$, $\xi>0$, extends to a continuous
  function on $[0,\infty)$ with $\rho(0)=\phi'(0)$, whence combining the positivity of $\rho$ on $(0,K]$ with that of $\phi'(0)$,
  as both being ensured by (\ref{phi}), we obtain that $\lambda(K):=\inf_{\xi\in [0,K]} \rho(\xi)$ is positive and satisfies
  the left inequality in (\ref{1.1}).
\qed
We next derive a space-time $L^2$ bound for $\ueps$, weigthed by the factor $\veps$ due to the lower bound from (\ref{1.1}),
by adapting a duality-based strategy which appears well-established in the analysis of semilinear parabolic problems,
but which also partially been pursued in some contexts of cross-diffusive systems related to (\ref{0}) 
(\cite{desvillettes}, \cite{taowin_M3AS}).
Unlike in most precedents, however, thanks to (\ref{uv}) the information thereby generated will here even include
corresponding integrability over the whole existence interval, thus implictly containing a certain decay information.
\begin{lem}\label{lem2}
  If (\ref{init}) holds, then there exists $C>0$ such that
  \be{2.2}
	\int_0^{\tme} \io \ueps^2 \veps \le C
	\qquad \mbox{for all } \eps\in (0,1).
  \ee
\end{lem}
\proof
  For $\vp\in L^1(\Om)$ we abbreviate $\ov{\vp}:=\frac{1}{|\Om|} \io \vp$, and
  we let $A$ denote the realization of $-\Del$ in $L^2_\perp(\Om):=\big\{\vp\in L^2(\Om) \ | \ \ov{\vp}=0\big\}$, with its domain
  given by $D(A):=\{\vp\in W^{2,2}(\Om) \cap L^2_\perp(\Om) \ | \ \frac{\pa\vp}{\pa\nu}=0 \mbox{ on } \pO\}$.
  Then $A$ is self-adjoint and positive, and an application of $A^{-1}$ to the identity
  \bas
	(\ueps-\ouz)_t &=&
	\Del \Big\{ \eps(\ueps-\ouz) + \big( \ueps\phi(\veps)-\ov{\ueps\phi(\veps)} \big) \Big\} \\
	&=& -A \Big\{ \eps(\ueps-\ouz) + \big( \ueps\phi(\veps)-\ov{\ueps\phi(\veps)} \big) \Big\},
	\qquad x\in\Om, \ t\in (0,\tme),
  \eas
  as implied by (\ref{0eps}) and (\ref{mass}), upon testing by $\ueps-\ouz$ shows that
  \bas
	\frac{1}{2} \frac{d}{dt} \io \Big| A^{-\frac{1}{2}} (\ueps-\ouz)\Big|^2
	&=& - \io \Big\{ \eps(\ueps-\ouz) + \big( \ueps\phi(\veps)-\ov{\ueps\phi(\veps)} \big) \Big\} \cdot (\ueps-\ouz) \\
	&=& -\eps\io (\ueps-\ouz)^2
	- \io \ueps^2 \phi(\veps)
	+ \ouz \io \ueps\phi(\veps)
	\qquad \mbox{for all } t\in (0,\tme),
  \eas
  because
  \bas
	\io \ov{\ueps\phi(\veps)} \cdot (\ueps-\ouz)=0
	\qquad \mbox{for all } t\in (0,\tme),
  \eas
  again due to (\ref{mass}).
  In view of Lemma \ref{lem1}, after integrating in time and dropping two favorably signed summands
  we thus obtain that with $K:=\|v_0\|_{L^\infty(\Om)}$ we have
  \bas
	\lambda(K) \int_0^t \io \ueps^2 \veps
	&\le& \frac{1}{2} \io \Big| A^{-\frac{1}{2}} (u_0-\ouz)\Big|^2
	+ \ouz \int_0^t \io \ueps \phi(\veps) \\
	&\le& \frac{1}{2} \io \Big| A^{-\frac{1}{2}} (u_0-\ouz)\Big|^2
	+ \Lambda(K) \ouz \int_0^t \io \ueps\veps
	\qquad \mbox{for all $t\in (0,\tme)$ and } \eps\in (0,1).
  \eas
  According to (\ref{uv}), this entails (\ref{2.2}) with an obvious choice of $C$.
\qed
This enables us to suitably control the interaction-driven contributions that appear in a standard first-oder testing
procedure applied to the second equation from (\ref{0eps}).
As we are assuming $\Om$ to be convex, corresponding boundary integrals are conveniently signed and hence the overall estimates
thereby gained again including the entire time range $(0,\tme)$.
An interesting question left open here is how far a large-time relaxation feature similar to that implicitly expressed 
in (\ref{3.2}) can be derived also
in more general domains; while our existence theory in the context of Theorem \ref{theo17} could readily be extended to such
settings by adaptations based on fairly well-established arguments, convexity seems more essential in the parts
in which Lemma \ref{lem3} will be applied in the large-time analysis addressing boundedness and stabilization
properties of solutions (cf.~Lemma \ref{lem5} and Lemma \ref{lem166}, for instance).
\begin{lem}\label{lem3}
  Assume (\ref{init}). Then there exists $C>0$ such that
  \be{3.2}
	\int_0^{\tme} \io \frac{|\na\veps|^4}{\veps^3} 	
	\le C
	\qquad \mbox{for all } \eps\in (0,1).
  \ee
\end{lem}
\proof
  By straightforward computation using (\ref{0eps}) and integration by parts (cf.~also \cite[Lemma 3.2]{win_CPDE2012}), we obtain
  the identity
  \bea{3.3}
	& & \hs{-30mm}
	\frac{1}{2} \frac{d}{dt} \io \frac{|\na\veps|^2}{\veps}
	+ \io \veps |D^2\ln\veps|^2
	+ \frac{1}{2} \io \frac{\ueps}{\veps}|\na\veps|^2 \nn\\
	&=& - \io \na\ueps\cdot\na\veps 
	+ \frac{1}{2} \io \frac{1}{\veps} \frac{\pa |\na\veps|^2}{\pa\nu}
	\qquad \mbox{for all } t\in (0,\tme),
  \eea
  where the rightmost summand is nonpositive, because $\frac{\pa |\na\veps|^2}{\pa\nu} \le 0$ on $\pO\times (0,\tme)$
  by convexity of $\Om$ (\cite{lions_ARMA}).
  As it is well-known (\cite[Lemma 3.3]{win_CPDE2012} \cite[Lemma 3.4]{ct_doubly_degenerate_smallsignal})
  that there exist positive constants $c_1$ and $c_2$ such that
  for all $\vp\in C^2(\bom)$ such that $\vp>0$ in $\bom$ and $\frac{\pa\vp}{\pa\nu}=0$ on $\pO$ we have
  \bas
	c_1 \io \frac{|\na\vp|^4}{\vp^3}
	\le \io \vp |D^2\ln\vp|^2
	\quad \mbox{and} \quad
	c_2 \io \frac{|D^2 \vp|^2}{\vp}
	\le \io \vp |D^2\ln\vp|^2,
  \eas
  from (\ref{3.3}) we thus infer that
  \be{3.4}
	4\frac{d}{dt} \io |\na\sqrt{\veps}|^2
	+ c_1 \io \frac{|\na\veps|^4}{\veps^3} + c_2 \io \frac{|D^2\veps|^2}{\veps}
	+ \io \frac{\ueps}{\veps} |\na\veps|^2
	\le - 2\io \na\ueps\cdot\na\veps
	\qquad \mbox{for all } t\in (0,\tme),
  \ee
  where now, after a further integration by parts, we may use Young's inequality to estimate
  \bas
	-2\io \na\ueps\cdot\na\veps
	&=& 2\io \ueps\Del\veps \\
	&\le& \frac{c_2}{n} \io \frac{|\Del\veps|^2}{\veps}
	+ \frac{n}{c_2} \io \ueps^2\veps \\
	&\le& c_2 \io \frac{|D^2\veps|^2}{\veps}
	+ \frac{n}{c_2} \io \ueps^2\veps
	\qquad \mbox{for all } t\in (0,\tme),
  \eas
  because $|\Del\veps|^2 \le n|D^2\veps|^2$. Therefore,
  \bas
	c_1 \int_0^t \io \frac{|\na\veps|^4}{\veps^3}
	+ \int_0^t \io \frac{\ueps}{\veps} |\na\veps|^2
	\le 4\io |\na \sqrt{v_0}|^2
	+ \frac{n}{c_2} \int_0^t \io \ueps^2\veps
	\qquad \mbox{for all } t\in (0,\tme),
  \eas
  so that (\ref{3.2}) results from Lemma \ref{lem2}.
\qed
To provide a prerequisite for a subsequent $L^p$ regularity argument concerning $\ueps$,
as the second of our key tools we now address the functional inequality announced in (\ref{fi}).
We underline that its derivation does actually not require any convexity hypothesis, but through the use
of a Sobolev embedding property it relies on the assumption that the spatial setting be one- or two-dimensional.
\begin{lem}\label{lem4}
  Let $n\le 2$ and $G\subset\R^n$ be a bounded domain with smooth boundary, and let 
  $p\ge 2$. Then there exists $C(p,G)>0$ such that for any $\vp\in C^1(\ov{G})$ and $\psi\in C^1(\ov{G})$ fulfilling
  $\vp\ge 0$ and $\psi>0$ in $\ov{G}$,
  \bea{4.1}
	\hs{-6mm}
	\ig \frac{\vp^p}{\psi} |\na\psi|^2
	&\le& \eta \ig \vp^{p-2} \psi |\na\vp|^2
	+ \eta \ig \vp\psi \nn\\
	& & + C(p,G) \cdot \Big(1+\frac{1}{\eta}\Big) \cdot \Bigg\{ \ig \vp^p + \bigg\{ \ig \vp \bigg\}^{2p-1} \Bigg\}
		\cdot \ig \frac{|\na\psi|^4}{\psi^3}
	\qquad \mbox{for all } \eta>0.
  \eea
\end{lem}
\proof
  Since we are assuming that $n\le 2$, and that thus $W^{1,1}(G)$ is continuously embedded into $L^2(G)$, 
  a corresponding Sobolev inequality yields $c_1(G)>0$ fulfilling
  \bas
	\|\rho\|_{L^2(G)} \le c_1(G) \|\na\rho\|_{L^1(G)} + c_1(G)\|\rho\|_{L^1(G)}
	\qquad \mbox{for all } \rho\in W^{1,1}(G),
  \eas
  so that since H\"older's and Young's inequalities imply that
  \bas
	c_1(G)\|\rho\|_{L^1(G)}
	&\le& c_1(G)\|\rho\|_{L^2(G)}^\frac{2p-2}{2p-1} \|\rho\|_{L^\frac{1}{p}(G)}^\frac{1}{2p-1} \\
	&=& \Big\{ \frac{1}{2} \|\rho\|_{L^2(G)} \Big\}^\frac{2p-2}{2p-1} 
		\cdot 2^\frac{2p-2}{2p-1} c_1(G) \|\rho\|_{L^\frac{1}{p}(G)}^\frac{1}{2p-1} \\
	&\le& \frac{1}{2} \|\rho\|_{L^2(G)}
	+ 2^{2p-2} c_1^{2p-1}(G) \|\rho\|_{L^\frac{1}{p}(G)}
	\qquad \mbox{for all } \rho\in L^2(G),
  \eas
  it follows that
  \bas
	\|\rho\|_{L^2(G)} \le c_2(p,G) \|\na\rho\|_{L^1(G)} + c_2(p,G)\|\rho\|_{L^\frac{1}{p}(G)}
	\qquad \mbox{for all } \rho\in W^{1,1}(G)
  \eas
  with $c_2(p,G):=\max\{2c_1(G) \, , \, (2c_1(G))^{2p-1}\}$.
  On the right-hand side of the estimate
  \be{4.2}
	\ig \frac{\vp^p}{\psi} |\na\psi|^2
	\le \bigg\{ \ig \frac{|\na\psi|^4}{\psi^3} \bigg\}^\frac{1}{2} \cdot
	\bigg\{ \ig \vp^{2p}\psi \bigg\}^\frac{1}{2},
  \ee
  valid whenever $0\le\vp\in C^1(\ov{G})$ and $0<\psi\in C^1(\ov{G})$ by the Cauchy-Schwarz inequality, we can therefore control
  the second factor according to
  \bea{4.3}
	\hs{-8mm}
	\bigg\{ \ig \vp^{2p}\psi \bigg\}^\frac{1}{2}
	&=& \|\vp^p \sqrt{\psi}\|_{L^2(\Om)} \nn\\
	&\le& 
	c_2(p,G)
	\ig \Big| p\vp^{p-1} \sqrt{\psi} \na\vp + \frac{\vp^p}{2\sqrt{\psi}} \na\psi \Big| 
	+ c_2(p,G)
	\cdot \bigg\{ \ig \vp \psi^\frac{1}{2p} \bigg\}^p \nn\\
	&\le& p c_2(p,G)
 	\ig \vp^{p-1} \sqrt{\psi} |\na\vp|
	+ \frac{c_2(p,G)}{2}
	\ig \frac{\vp^p}{\sqrt{\psi}} |\na\psi| 
	+ c_2(p,G)
	\cdot \bigg\{ \ig \vp \psi^\frac{1}{2p} \bigg\}^p.
  \eea
  Here three applications of the H\"older inequality show that
  \bas
	p c_2(p,G)
	\ig \vp^{p-1} \sqrt{\psi} |\na\vp|
	\le pc_2(p,G)
	\cdot \bigg\{ \ig \vp^p \bigg\}^\frac{1}{2} \cdot
	\bigg\{ \ig \vp^{p-2} \psi |\na\vp|^2 \bigg\}^\frac{1}{2}
  \eas
  and
  \bas
	\frac{c_2(p,G)}{2} 
	\ig \frac{\vp^p}{\sqrt{\psi}} |\na\psi| 
	\le 
	\frac{c_2(p,G)}{2}
	\cdot \bigg\{ \ig \vp^p \bigg\}^\frac{1}{2} \cdot 
	\bigg\{ \ig \frac{\vp^p}{\psi} |\na\psi|^2 \bigg\}^\frac{1}{2}
  \eas
  as well as
  \bas	
	c_2(p,G)
	\cdot \bigg\{ \ig \vp \psi^\frac{1}{2p} \bigg\}^p
	&=& 
	c_2(p,G)
	\cdot \bigg\{ \ig (\vp\psi)^\frac{1}{2p} \cdot \vp^\frac{2p-1}{2p} \bigg\}^p \\
	&\le& 
	c_2(p,G)
	\cdot \bigg\{ \ig \vp\bigg\}^\frac{2p-1}{2} \cdot \bigg\{ \ig \vp\psi\bigg\}^\frac{1}{2}.
  \eas
  Inserting (\ref{4.3}) into (\ref{4.2}) and using Young's inequality, we hence infer that for each $\eta>0$,
  \bas
	\ig \frac{\vp^p}{\psi} |\na\psi|^2
	&\le& p
	c_2(p,G) 
	\cdot \bigg\{ \ig \frac{|\na\psi|^4}{\psi^3} \bigg\}^\frac{1}{2} \cdot
	\bigg\{ \ig \vp^p\bigg\}^\frac{1}{2} \cdot
	\bigg\{ \ig \vp^{p-2} \psi |\na\vp|^2 \bigg\}^\frac{1}{2} \\
	& & + 
	\frac{c_2(p,G)}{2} 
	\cdot \bigg\{ \ig \frac{|\na\psi|^4}{\psi^3} \bigg\}^\frac{1}{2} \cdot
	\bigg\{ \ig \vp^p \bigg\}^\frac{1}{2} \cdot
	\bigg\{ \ig \frac{\vp^p}{\psi} |\na\psi|^2 \bigg\}^\frac{1}{2} \\
	& & + 
	c_2(p,G) 
	\cdot \bigg\{ \ig \frac{|\na\psi|^4}{\psi^3} \bigg\}^\frac{1}{2} \cdot 
	\bigg\{ \ig \vp \bigg\}^\frac{2p-1}{2} \cdot
	\bigg\{ \ig \vp\psi\bigg\}^\frac{1}{2} \\
	&\le& \frac{\eta}{2} \ig \vp^{p-2} \psi |\na\vp|^2
	+ \frac{p^2 c_2^2(p,G)}{2\eta} 
	\cdot \bigg\{ \ig \vp^p\bigg\} \cdot \ig \frac{|\na\psi|^4}{\psi^3} \\
	& & + \frac{1}{2} \ig \frac{\vp^p}{\psi} |\na\psi|^2
	+ \frac{c_2^2(p,G)}{8} 
	\cdot \bigg\{ \ig \vp^p \bigg\} \cdot \ig \frac{|\na\psi|^4}{\psi^3} \\
	& & + \frac{\eta}{2} \ig \vp\psi
	+ \frac{c_2^2(p,G)}{2\eta} 
	\cdot \bigg\{ \ig \vp\bigg\}^{2p-1} \cdot \ig \frac{|\na\psi|^4}{\psi^3},
  \eas
  which readily implies (\ref{4.1}) with $C(p,G):=p^2 c_2^2(p,G)$.
\qed
We are now prepared to make sure that despite the diffusion degeneracy in the first equation from (\ref{0eps}),
the respective first solution components remain bounded with respect to the norm in any $L^p$ space with $p\ge 2$.
Besides on Lemma \ref{lem4}, our derivation of this will make substantial use of the corresponding decay features
expressed in (\ref{3.2}) and, again, in (\ref{uv}).
\begin{lem}\label{lem5}
  Given any $p\ge 2$, one can pick $C(p)>0$ such that if (\ref{init}) holds, then
  \be{5.1}
	\io \ueps^p(\cdot,t) \le C(p)
	\qquad \mbox{for all $t\in (0,\tme)$ and } \eps\in (0,1).
  \ee
\end{lem}
\proof
  We first employ Lemma \ref{lem3} to fix $c_1>0$ such that
  \be{5.2}
	\int_0^{\tme} \io \frac{|\na\veps|^4}{\veps^3} \le c_1
	\qquad \mbox{for all } \eps\in (0,1),
  \ee
  and to make adequate use of this together with the outcome of Lemma \ref{lem4},
  we utilize Young's inequality when testing the first equation in (\ref{0eps}) by $\ueps^{p-1}$
  to find that
  \bea{5.22}
	\frac{d}{dt} \io \ueps^p
	&=& p\io \ueps^{p-1} \Del \big\{ \eps\ueps + \ueps\phi(\veps) \big\} \nn\\
	&=& -p(p-1) \eps \io \ueps^{p-2} |\na\ueps|^2
	- p(p-1) \io \ueps^{p-2} \phi(\veps) |\na\ueps|^2  \nn\\
	& & - p(p-1) \io \ueps^{p-1} \phi'(\veps) \na\ueps\cdot\na\veps \nn\\
	&\le& -\frac{p(p-1)}{2} \io \ueps^{p-2} \phi(\veps) |\na\ueps|^2 \nn\\
	& & + \frac{p(p-1)}{2} \io \ueps^p \frac{\phi'^2(\veps)}{\phi(\veps)} |\na\veps|^2
	\quad \mbox{for all } t\in (0,\tme).
  \eea
  Here, abbreviating $K:=\|v_0\|_{L^\infty(\Om)}$ we may draw on Lemma \ref{lem1} in estimating
  \bas
	\phi(\veps)\ge \lambda(K)\veps
	\quad \mbox{and} \quad
	\frac{\phi'^2(\veps)}{\phi(\veps)} \le \frac{\Lambda^2(K)}{\lambda(K) \veps}
	\qquad \mbox{in } \Om\times (0,\tme),
  \eas
  so that since an application of Lemma \ref{lem4} to $\eta:=\min\{ \frac{p(p-1)\Lambda^2(K)}{2\lambda(K)} \, , \, 1\}$
  provides $c_2(p)>0$ such that for all $\vp\in C^1(\bom)$ and $\psi\in C^1(\bom)$ with
  $\vp\ge 0$ and $\psi>0$ in $\bom$ we have
  \bas
	\frac{p(p-1)\Lambda^2(K)}{2\lambda(K)} \io \frac{\vp^p}{\psi} |\na\psi|^2
	&\le& \frac{p(p-1) \lambda(K)}{2} \io \vp^{p-2} \psi |\na\vp|^2
	+ \io \vp\psi \\
	& & + c_2(p) \cdot \Bigg\{ \io \vp^p + \bigg\{ \io \vp\bigg\}^{2p-1} \Bigg\} \cdot \io \frac{|\na\psi|^4}{\psi^3},
  \eas
  thanks to (\ref{mass}) this entails that for all $t\in (0,\tme)$ and $\eps\in (0,1)$,
  \bas
	\frac{d}{dt} \io \ueps^p
	\le \io \ueps\veps
	+ c_2(p) \cdot \Bigg\{ \io \ueps^p + \bigg\{ \io u_0 \bigg\}^{2p-1} \Bigg\} \cdot \io \frac{|\na\veps|^4}{\veps^3}.
  \eas
  For each $\eps\in (0,1)$, the functions given by
  \bas
	\yeps(t):= \io \ueps^p(\cdot,t) + \bigg\{ \io u_0 \bigg\}^{2p-1},
	\qquad t\in [0,\tme),
  \eas
  as well as
  \bas
	g_\eps(t):=\io \ueps(\cdot,t)\veps(\cdot,t)
	\quad \mbox{and} \quad
	h_\eps(t):=c_2(p) \io \frac{|\na\veps(\cdot,t)|^4}{\veps^3(\cdot,t)},
	\qquad t\in (0,\tme),
  \eas
  thus satisfy
  \bas
	\yeps'(t) \le g_\eps(t) + h_\eps(t) \yeps(t)
	\qquad \mbox{for all } t\in (0,\tme),
  \eas
  which upon an ODE comparison argument implies that
  \be{5.3}
	\yeps(t) \le \yeps(0) e^{\int_0^t h_\eps(s) ds} 
	+ \int_0^t e^{\int_s^t h_\eps(\sigma)d\sigma} g_\eps(s) ds 
	\qquad \mbox{for all } t\in (0,\tme).
  \ee
  Since
  \bas
	\int_s^t h_\eps(\sigma) d\sigma
	\le c_1 c_2(p)
	\qquad \mbox{for all $t\in (0,\tme)$, $s\in [0,t)$ and } \eps\in (0,1)
  \eas
  by (\ref{5.2}), and since
  \bas
	\int_0^t g_\eps(s) ds \le \io v_0
	\qquad \mbox{for all $t\in (0,\tme)$ and } \eps\in (0,1)
  \eas
  due to (\ref{uv}), from (\ref{5.3}) we thus obtain that
  \bas
	\io \ueps^p(\cdot,t) 
	\le \Bigg\{ \io u_0^p + \bigg\{ \io u_0\bigg\}^{2p-1} \Bigg\} \cdot e^{c_1 c_2(p)} + e^{c_1 c_2(p)} \io v_0
	\qquad \mbox{for all $t\in (0,\tme)$ and } \eps\in (0,1)
  \eas
  to conclude as intended.
\qed
\mysection{Analysis of (\ref{0}): Positivity properties of $\veps$ and higher order estimates}
Thanks to Lemma \ref{lem5}, we are now in the position to draw the intended conclusion from Proposition \ref{prop6},
and to thereby obtain a pointwise lower estimate for the second solution components which indeed is uniform with respect
to the approximation parameter:
\begin{cor}\label{cor7}
  Assume (\ref{init}). 
  Then for all $T>0$ and $\tau\in (0,T)$ there exists $C(T,\tau)>0$ such that
  \be{7.1}
	\veps(x,t) \ge C(T,\tau)
	\qquad \mbox{for all $x\in\Om, \ t\in (\tau,T)\cap (0,\tme)$ and } \eps\in (0,1).
  \ee
\end{cor}
\proof
  Since $v_0\not\equiv 0$, this immediately results from Proposition \ref{prop6} upon applying Lemma \ref{lem5} to, e.g., $p:=2$.
\qed
Independently from the latter, through standard parabolic regularity arguments the outcome of Lemma \ref{lem5} furthermore
entails uniform bounds for the taxis gradients in (\ref{0eps}):
\begin{lem}\label{lem8}
  If (\ref{init}) holds, then there exists $C>0$ such that
  \be{8.1}
	\|\na\veps(\cdot,t)\|_{L^\infty(\Om)} \le C
	\qquad \mbox{for all $t\in (0,\tme)$ and } \eps\in (0,1).
  \ee
\end{lem}
\proof
  According to well-known smoothing properties of the Neumann heat semigroup $(e^{t\Del})_{t\ge 0}$ on $\Om$ (\cite{win_JDE2010}),
  fixing any $p>2$ we can find $c_1>0$ such that for all $t\in  (0,\tme)$ and $\eps\in (0,1)$,
  \bea{8.3}
	\|\na\veps(\cdot,t)\|_{L^\infty(\Om)}
	&=& \bigg\| \na e^{t(\Del-1)} v_0
	- \int_0^t \na e^{(t-s)(\Del-1)} \Big\{ \ueps(\cdot,s)\veps(\cdot,s)-\veps(\cdot,s)\Big\} ds \bigg\|_{L^\infty(\Om)} \nn\\
	&\le& c_1 \|v_0\|_{W^{1,\infty}(\Om)} \nn\\[2mm]
	& & + c_1\int_0^t \Big(1+(t-s)^{-\frac{1}{2}-\frac{n}{2p}}\Big) e^{-(t-s)} 
		\big\| \ueps(\cdot,s)\veps(\cdot,s)-\veps(\cdot,s)\big\|_{L^p(\Om)} ds.
  \eea
  Since (\ref{vinfty}) implies that
  \bas
	& & \hs{-20mm}
	\big\| \ueps(\cdot,s)\veps(\cdot,s)-\veps(\cdot,s)\big\|_{L^p(\Om)} \\
	&\le& \|\ueps(\cdot,s)\|_{L^p(\Om)} \|\veps(\cdot,s)\|_{L^\infty(\Om)} + |\Om|^\frac{1}{p} \|\veps(\cdot,s)\|_{L^\infty(\Om)}
		\\
	&\le& c_2 \|v_0\|_{L^\infty(\Om)} + |\Om|^\frac{1}{p} \|v_0\|_{L^\infty(\Om)}
	\qquad \mbox{for all $s\in (0,\tme)$ and } \eps\in (0,1),
  \eas
  with $c_2:=\sup_{\eps\in (0,1)} \sup_{t\in (0,\tme)} \|\ueps(\cdot,s)\|_{L^p(\Om)}$ being finite by Lemma \ref{lem5},
  from (\ref{8.3}) we directly obtain (\ref{8.1}).
\qed
Relying on information on actual non-degeneracy of diffusion in (\ref{0eps}), as implied by Corollary \ref{cor7} 
throughout any region of the form $\Om\times \big((\tau,T) \cap (0,\tme)\big)$ with $0<\tau<T$,
by means of a straightforward temporal cut-off procedure 
we can now utilize Lemma \ref{lem8} to establish local-in-time $L^\infty$ bounds for $\ueps$ through the outcome
of a Moser-type iterative reasoning.
\begin{lem}\label{lem9}
  Suppose that (\ref{init}) holds. Then for all $T>0$ and $\tau\in (0,T)$ one can find $C(T,\tau)>0$ such that
  \be{9.1}
	\|\ueps(\cdot,t)\|_{L^\infty(\Om)} \le C(T,\tau)
	\qquad \mbox{for all $t\in (\tau,T) \cap (0,\tme)$ and } \eps\in (0,1).
  \ee
\end{lem}
\proof
  We fix $\zeta\in C^\infty([0,\infty))$ such that $\zeta\equiv 0$ on $[0,\frac{\tau}{2}]$ and $\zeta\equiv 1$ on $[\tau,\infty)$,
  and then from (\ref{0eps}) we obtain that $\weps(x,t):=\zeta(t)\cdot\ueps(x,t)$, $(x,t)\in\bom\times [0,\tme)$, $\eps\in (0,1)$,
  satisfies
  \be{9.2}
	w_{\eps t} = \na \cdot \big( D_\eps(x,t)\na\weps\big) + \na\cdot f_\eps(x,t) + g_\eps(x,t),
	\qquad x\in \Om, \, t\in (0,\tme), \ \eps\in (0,1),
  \ee
  where 
  \bas
	D_\eps(x,t):=\eps+\phi\big(\veps(x,t)\big),
	\quad
	f_\eps(x,t):=\zeta(t)\ueps(x,t) \phi'\big(\veps(x,t)\big) \na\veps(x,t)
	\quad \mbox{and} \quad
	g_\eps(x,t):=\zeta'(t) \ueps(x,t)
  \eas
  for $(x,t)\in\Om\times (0,\tme)$ and $\eps\in (0,1)$.
  Here since
  \bas
	\lambda(\|v_0\|_{L^\infty(\Om)}) \cdot \veps \le D_\eps(x,t) \le 1+\Lambda(\|v_0\|_{L^\infty(\Om)}) \cdot 
		\|v_0\|_{L^\infty(\Om)}
	\ \mbox{for all $x\in \Om, t\in (0,\tme)$ and } \eps\in (0,1)
  \eas
  by Lemma \ref{lem1} and (\ref{vinfty}), using Corollary \ref{cor7} we infer the existence of $c_1(T,\tau)>0$ and $c_2>0$ such that
  \bas
	c_1(T,\tau) \le D_\eps(x,t) \le c_2
	\qquad \mbox{for all $x\in \Om, t\in (\frac{\tau}{2},T) \cap (0,\tme)$ and } \eps\in (0,1).  
  \eas
  Since, apart from that, a combination of Lemma \ref{lem5} with (\ref{vinfty}) and Lemma \ref{lem8} shows that
  \bas
	\sup_{\eps\in (0,1)} \sup_{t\in (0,\tme)} \Big\{ \|\weps(\cdot,t)\|_{L^p(\Om)} + \|f_\eps(\cdot,t)\|_{L^p(\Om)}
		+ \|g_\eps(\cdot,t)\|_{L^p(\Om)} \Big\} <\infty
	\qquad \mbox{for all } p\in [2,\infty),
  \eas
  and since $\weps(\cdot,\frac{\tau}{2})\equiv 0$ for all $\eps\in (0,1)$ according to our choice of $\zeta$,
  an application of Lemma A.1 in \cite{taowin_subcrit} yields $c_3(T,\tau)>0$ such that
  \bas
	\|\weps(\cdot,t)\|_{L^\infty(\Om)}
	\le c_3(T,\tau)
	\qquad \mbox{for all $t\in (\frac{\tau}{2},T)\cap (0,\tme)$ and } \eps\in (0,1).
  \eas
  As $\weps(\cdot,t) \equiv \ueps(\cdot,t)$ in $\Om$ for all $t\in (\tau,T)\cap (0,\tme)$ and $\eps\in (0,1)$, this implies (\ref{9.1}).
\qed
The latter especially rules out any blow-up in the approximate problems:
\begin{lem}\label{lem10}
  If (\ref{init}) holds, then $\tme=+\infty$ for all $\eps\in (0,1)$.
\end{lem}
\proof
  This immediately follows from Lemma \ref{lem9} when combined with (\ref{ext}).
\qed
Apart from that, Lemma \ref{lem9} can be combined with Lemma \ref{lem8} in the course of an essentially straightforward
bootstrap procedure so as to yield temporally local higher-order regularity properties.
\begin{lem}\label{lem11}
  Assume (\ref{init}). Then for all $T>0$ and any $\tau\in (0,T)$ there exist $\theta=\theta(T,\tau)\in (0,1)$ and
  $C(T,\tau)>0$ such that
  \be{11.1}
	\|\ueps\|_{C^{2+\theta,1+\frac{\theta}{2}}(\bom\times [\tau,T])} \le C(T,\tau)
	\qquad \mbox{for all } \eps\in (0,1)
  \ee
  and
  \be{11.2}
	\|\veps\|_{C^{2+\theta,1+\frac{\theta}{2}}(\bom\times [\tau,T])} \le C(T,\tau)
	\qquad \mbox{for all } \eps\in (0,1).
  \ee
\end{lem}
\proof
  We rewrite the first equation from (\ref{0eps}) according to
  \bas
	u_{\eps t} = \na \cdot A_\eps(x,t,\na\ueps),
	\qquad x\in\Om, \ t>0, \ \eps\in (0,1),
  \eas
  with
  \bas
	A_\eps(x,t,\xi):=\eps\xi + \phi\big(\veps(x,t)\big) \xi
	+ \phi'\big(\veps(x,t)\big) \ueps(x,t)\na\veps(x,t),
	\qquad (x,t,\xi)\in\Om\times (0,\infty)\times\R, \ \eps\in (0,1),
  \eas
  and employ Corollary \ref{cor7} along with Lemma \ref{lem1}, (\ref{vinfty}), Lemma \ref{lem9} and Lemma \ref{lem8} to find
  $c_1(T,\tau)>0$, $c_2(T,\tau)>0$ and $c_3(T,\tau)>0$ such that whenever $\eps\in (0,1)$,
  \bas
	A_\eps(x,t,\xi) \cdot \xi \ge c_1(T,\tau) |\xi|^2 - c_2(T,\tau)
	\qquad \mbox{for all $(x,t,\xi)\in\Om\times (\frac{\tau}{8},T)\times\R^n$ and } \eps\in (0,1)
  \eas
  and
  \bas
	|A_\eps(,x,t,\xi)|
	\le c_3(T,\tau) |\xi| + c_3(T,\tau)
	\qquad \mbox{for all $(x,t,\xi)\in\Om\times (\frac{\tau}{8},T)\times\R^n$ and } \eps\in (0,1).
  \eas
  Again based on Lemma \ref{lem9}, by means of a standard result on H\"older regularity of bounded solutions to scalar parabolic
  equations (\cite{porzio_vespri}) we thus obtain $\theta_1=\theta_1(T,\tau)\in (0,1)$ and $c_4(T,\tau)>0$ such that
  \bas
	\|\ueps\|_{C^{\theta_1,\frac{\theta_1}{2}}(\bom\times [\frac{\tau}{4},T])} \le c_4(T,\tau)
	\qquad \mbox{for all } \eps\in (0,1),
  \eas
  whereupon parabolic Schauder theory applies to the second equation from (\ref{0eps}) so as to yield 
  $\theta_2=\theta_2(T,\tau)\in (0,1)$ and $c_5(T,\tau)>0$ fulfilling
  \be{11.3}
	\|\veps\|_{C^{2+\theta_2,1+\frac{\theta_2}{2}}(\bom\times [\frac{\tau}{2},T])} \le c_5(T,\tau)
	\qquad \mbox{for all } \eps\in (0,1).
  \ee
  This information in turn enables us to go back to the first equation in (\ref{0eps}), now written in the form
  \bas
	u_{\eps t} = \big\{ \eps+\phi(\veps)\big\} \Delta \ueps
	+ \big\{ 2\phi'(\veps)\na\veps\big\}\cdot \na\ueps
	+ \big\{ \phi'(\veps) \Del\veps + \phi''(\veps)|\na\veps|^2 \big\} \ueps,
	\qquad x\in\Om, \, t>0, \, \eps\in (0,1),
  \eas
  to conclude again from parabolic Schauer theory and the stimates provided by Corollary \ref{cor7}, (\ref{vinfty})
  and Lemma \ref{lem8} that (\ref{11.1}) holds with some $\theta=\theta(T,\tau)\in (0,1)$ and $C(T,\tau)>0$.
  In view of (\ref{11.3}), the proof thereby becomes complete.
\qed
As a last preparation for our limit passage, let us once more go back to Lemma \ref{lem5} to obtain the following 
information on H\"older regularity of $\veps$ down to the temporal origin.
\begin{lem}\label{lem12}
  If (\ref{init}) is satisfied, then for each $T>0$ there exist $\theta=\theta(T)\in (0,1)$ and
  $C(T)>0$ such that
  \be{12.1}
	\|\veps\|_{C^{\theta,\frac{\theta}{2}}(\bom\times [0,T])} \le C(T)
	\qquad \mbox{for all } \eps\in (0,1).
  \ee
\end{lem}
\proof
  This immediately follows from standard parabolic regularity theory (\cite{porzio_vespri}) after applying Lemma \ref{lem5}
  to any fixed $p\ge 2$.
\qed
A solution of (\ref{0}) in the flavor of the statement from Theorem \ref{theo17} can now be obtained by a standard
extraction process, followed by a suitably arranged argument asserting continuity of the corresponding first component
with respect to weak $L^p$ topologies.	
\begin{lem}\label{lem13}
  Assume (\ref{init}). Then there exist $(\eps_j)_{j\in\N} \subset (0,1)$ as well as functions $u$ and $v$ on $\bom\times (0,\infty)$
  such that $\eps_j\searrow 0$ as $j\to\infty$, that (\ref{17.1}) holds with $u> 0$ and $v>0$ in $\bom\times (0,\infty)$, and that
  \begin{eqnarray}
	& & \ueps\to u
	\qquad \mbox{in } C^{2,1}_{loc}(\bom\times (0,\infty)), 
	\label{13.2} \\
	& & \ueps\wto u
	\qquad \mbox{in } L^p_{loc}(\bom\times [0,\infty)) \qquad \mbox{for all } p\ge 1, 
	\label{13.3} \\
	& & \veps\to v
	\qquad \mbox{in } C^0_{loc}(\bom\times [0,\infty)) \mbox{ and in } C^{2,1}_{loc}(\bom\times (0,\infty)), 
	\qquad \mbox{and that}
	\label{13.4} \\
	& & \na\veps \wsto \na v
	\qquad \mbox{in } L^\infty(\Om\times (0,\infty))
	\label{13.5}
  \end{eqnarray}
  as $\eps=\eps_j\searrow 0$.
  In the classical sense, these functions satisfy $u_t=\Del (u\phi(v))$ and $v_t=\Del v-uv$ in $\Om\times (0,\infty)$
  with $\frac{\pa u}{\pa\nu}=\frac{\pa v}{\pa \nu}=0$ on $\pO\times (0,\infty)$ and $v(x,0)=v_0(x)$ for all $x\in\Om$,
  and moreover (\ref{17.2}) holds.
  Apart from that,
  \be{13.66}
	\io u(\cdot,t)=\io u_0
	\quad \mbox{as well as} \quad
	\|v(\cdot,t)\|_{L^\infty(\Om)} \le \|v_0\|_{L^\infty(\Om)}
	\qquad \mbox{for all } t>0,
  \ee
  and
  \be{13.67}
	\int_0^\infty \io uv \le \io v_0.
  \ee
\end{lem}
\proof
  The existence of $(\eps_j)_{j\in\N}$ and nonnegative functions $u$ and $v$ 
  with the properties in (\ref{17.1}) and (\ref{13.2})-(\ref{13.5})
  follows from Lemma \ref{lem11}, Lemma \ref{lem12}, Lemma \ref{lem5} and Lemma \ref{lem8} by means of a straightforward extraction
  procedure, 
  whereupon the claimed classical solution features can then immediately be verified by taking $\eps=\eps_j \searrow 0$ in (\ref{0eps})
  and using (\ref{13.2}), (\ref{13.4}) and the continuity of $\phi, \phi'$ and $\phi''$.
  Strict positivity of $u$ and $v$ throughout $\bom\times (0,\infty)$ can then a posteriori be deduced by applying 
  the classical strong maximum principle to the identities $v_t=\Del v-uv$ and $u_t=\Del \big(u\phi(v)\big)$,
  while (\ref{13.66}) and (\ref{13.67}) result from (\ref{mass}), (\ref{vinfty}) and
  (\ref{uv}) in conjunction with (\ref{13.2}), (\ref{13.4}) and Fatou's lemma.\abs
  It thus remains to derive
  the initial trace feature expressed in (\ref{17.2}) for each $p\ge 1$, and
  to achieve this, assuming without loss of generality that $p>1$ we let $\psi\in (L^p(\Om))^\star \cong L^\frac{p}{p-1}(\Om)$
  and $\eta>0$ be given and pick any $\psi_\eta\in C_0^\infty(\Om)$ such that, in accordance with Lemma \ref{lem5} and (\ref{13.2}),
  we have
  \be{13.7}
	\|u(\cdot,t)\|_{L^p(\Om)} \cdot \|\psi-\psi_\eta\|_{L^\frac{p}{p-1}(\Om)}
	\le \frac{\eta}{3}
	\quad \mbox{for all } t>0
	\qquad \mbox{and} \qquad
	\|u_0\|_{L^p(\Om)} \cdot \|\psi-\psi_\eta\|_{L^\frac{p}{p-1}(\Om)}
	\le \frac{\eta}{3}.
  \ee
  We thereafter choose $t_\eta\in (0,1)$ suitably small such that with $\Lambda$ taken from Lemma \ref{lem1} we have
  \be{13.8}
	\bigg\{ \io u_0 \bigg\} \cdot \Lambda (\|v_0\|_{L^\infty(\Om)}) \cdot \|v_0\|_{L^\infty(\Om)} \cdot 
		\|\Del\psi_\eta\|_{L^\infty(\Om)} \cdot t_\eta \le \frac{\eta}{3},
  \ee
  and we claim that these selections guarantee that
  \be{13.87}
	\bigg| \io u(\cdot,t)\psi - \io u_0 \psi \bigg| \le \eta
	\qquad \mbox{for all } t\in (0,t_\eta).
  \ee
  In fact, since $\psi_\eta$ belongs to $C_0^\infty(\Om)$, when testing the first equation in (\ref{0eps}) against $\psi_\eta$
  we do not encounter nontrivial boundary integrals and hence obtain that
  \be{13.9}
	\io \ueps(\cdot,t)\psi_\eta - \io u_0\psi_\eta
	= \eps\int_0^t \io \ueps\Del\psi_\eta
	+ \int_0^t\io \ueps\phi(\veps)\Del\psi_\eta
	\qquad \mbox{for all $t>0$ and } \eps\in (0,1).
  \ee
  Here, fixing any $t\in (0,t_\eta)$ we may invoke (\ref{13.2}) to see that
  \bas
	\io \ueps(\cdot,t)\psi_\eta \to \io u(\cdot,t)\psi_\eta
	\qquad \mbox{as } \eps=\eps_j\searrow 0,
  \eas
  while combining (\ref{13.3}) with (\ref{13.4}) and the continuity of $\phi$ readily implies that
  \bas
	\eps\int_0^t \io \ueps\Del\psi_\eta \to 0
	\quad \mbox{and} \quad
	\int_0^t\io \ueps\phi(\veps)\Del\psi_\eta
	\to \int_0^t \io u\phi(v) \Del\psi_\eta
	\qquad \mbox{as } \eps=\eps_j\searrow 0.
  \eas
  Accordingly, (\ref{13.9}) entails that due to (\ref{13.66}), Lemma \ref{lem1} and (\ref{13.8}),
  \bas
	& & \hs{-20mm}
	\bigg| \io u(\cdot,t)\psi_\eta - \io u_0 \psi_\eta \bigg| \\
	&=& \bigg| \int_0^t \io u\phi(v)\Del\psi_\eta \bigg| \\
	&\le& \int_0^t \|u(\cdot,s)\|_{L^1(\Om)} \|\phi(v(\cdot,s))\|_{L^\infty(\Om)} \|\Del\psi_\eta\|_{L^\infty(\Om)} ds \\
	&\le& \bigg\{ \io u_0 \bigg\} \cdot \Lambda(\|v_0\|_{L^\infty(\Om)}) \cdot \|v_0\|_{L^\infty(\Om)} \cdot
	\|\Del\psi_\eta\|_{L^\infty(\Om)} \cdot t \\[2mm]
	&\le& \frac{\eta}{3},
  \eas
  because $t\in (0,t_\eta)$. 
  In view of (\ref{13.7}), we thus obtain that, indeed,
  \bas
	& & \hs{-20mm}
	\bigg| \io u(\cdot,t)\psi - \io u_0 \psi \bigg| \\
	&=& \Bigg|
	\io u(\cdot,t) \cdot (\psi-\psi_\eta)
	+ \bigg\{ \io u(\cdot,t) \psi_\eta
	- \io u_0 \psi_\eta \bigg\}
	+ \io u_0 \cdot (\psi_\eta-\psi) \Bigg| \\[1mm]
	&\le& \frac{\eta}{3} + \frac{\eta}{3} + \frac{\eta}{3}=\eta,
  \eas
  and that hence the verification of (\ref{13.87}), as thereby achieved, completes the proof.
\qed
Our main result concerning global solvability in (\ref{0}) has thereby been accomplished:\abs
\proofc of Theorem \ref{theo17}. \quad
  We only need to take $(u,v)$ as provided by Lemma \ref{lem13}.
\qed
\mysection{Large time behavior in (\ref{0}). Proof of Theorems \ref{theo18} and \ref{theo19}}\label{sect_largetime}
Our analysis of the large time behavior in (\ref{0}) is rooted in the following consequence of (\ref{13.67})
on the total variation of $u$ when considered as a $W_N^{2,\infty}(\Om)$-valued function over $[0,\infty)$.
Here and below, for definiteness in our corresponding argument we shall let the
Banach space 
$W_N^{2,\infty}(\Om)$, as introduced before Theorem \ref{theo18},
be equipped with the norm given by $\|\vp\|_{W^{2,\infty}(\Om)}:=\max_{|\alpha|\le 2} \|D^\alpha\vp\|_{L^\infty(\Om)}$,
$\vp\in W_N^{2,\infty}(\Om)$.
\begin{lem}\label{lem14}
  Let $K>0$. Then there exists $C(K)>0$ with the property that if (\ref{init}) holds with $\|v_0\|_{L^\infty(\Om)} \le K$,  
  for any choice of $(t_k)_{k\in\N}\subset [0,\infty)$ such that $t_{k+1} \ge t_k$ for all $k\in\N$, we have
  \be{14.1}
	\sum_{k\in\N} \big\| u(\cdot,t_{k+1})-u(\cdot,t_k)\big\|_{(W_N^{2,\infty}(\Om))^\star} \le C(K) \io v_0.
  \ee
\end{lem}
\proof
  For fixed $\psi\in W_N^{2,\infty}(\Om)$, an integration by parts in (\ref{0eps}) shows that
  \bas
	& & \hs{-20mm}
	\io \ueps(\cdot,t_{k+1}) \cdot \psi - \io \ueps(\cdot,t_k) \cdot\psi \\
	&=& \eps \int_{t_k}^{t_{k+1}} \io \ueps\Del\psi 
	+ \int_{t_k}^{t_{k+1}} \io \ueps\phi(\veps)\Del\psi
	\qquad \mbox{for all $k\in\N$ and } \eps\in (0,1),
  \eas
  and that hence, by (\ref{13.2}), (\ref{13.3}), (\ref{13.4}) and the continuity of $\phi$,
  \bas
	\io u(\cdot,t_{k+1})\cdot\psi - \io u(\cdot,t_k)\cdot\psi
	= \int_{t_k}^{t_{k+1}} \io u\phi(v)\Del\psi
	\qquad \mbox{for all } k\in\N.
  \eas
  Since $\phi(v)\le \Lambda(K) v$ according to Lemma \ref{lem1}, (\ref{13.4}) and our assumption, this implies that
  \bas
	\bigg| \io \Big\{ u(\cdot,t_{k+1})-u(\cdot,t_k) \Big\} \cdot \psi \bigg|
	\le \Lambda(K) \|\Del\psi\|_{L^\infty(\Om)} \int_{t_k}^{t_{k+1}} \io uv
	\qquad \mbox{for all } k\in\N,
  \eas
  so that estimating $\|\Del\psi\|_{L^\infty(\Om)} \le n\|\psi\|_{W^{2,\infty}(\Om)}$ we obtain that
  \bas
	\big\| u(\cdot,t_{k+1})-u(\cdot,t_k)\big\|_{(W_N^{2,\infty}(\Om))^\star} 
	\le n\Lambda(K) \int_{t_k}^{t_{k+1}} \io uv
	\qquad \mbox{for all } k\in\N
  \eas
  and thus
  \bas
	\sum_{k\in\N} \big\| u(\cdot,t_{k+1})-u(\cdot,t_k)\big\|_{(W_N^{2,\infty}(\Om))^\star} 
	\le n\Lambda(K) \int_0^\infty \io uv,
  \eas
  because $(t_k,t_{k+1})\cap (t_l,t_{l+1})=\emptyset$ for all $k\in\N$ and $l\in\N$ with $k\ne l$.
  The claim therefore results upon recalling (\ref{13.67}).
\qed
Thanks to the quantitative dependence on $v_0$, this does not only imply large time stabilization of each individual
trajectory in its first component, but it moreover provides some information on the distance between the associated limit 
and the initial data.
\begin{lem}\label{lem15}
  Let $K>0$. Then there exists $\Xi(K)>0$ such that whenever (\ref{init}) holds with $\|v_0\|_{L^\infty(\Om)} \le K$,  
  the function $u$ obtained in Lemma \ref{lem13} has the property that
  \be{15.1}
	u(\cdot,t) \to u_\infty
	\quad \mbox{in }
	\big(W_N^{2,\infty}(\Om)\big)^\star
	\qquad \mbox{as } t\to\infty,
  \ee
  with some $u_\infty\in \big(W_N^{2,\infty}(\Om)\big)^\star$ which satisfies
  \be{15.2}
	\|u_\infty-u_0\|_{(W_N^{2,\infty}(\Om))^\star} \le \Xi(K) \io v_0.
  \ee
\end{lem}
\proof
  Given any unbounded $(t_k)_{k\in\N} \subset (0,\infty)$ such that $t_{k+1}>t_k$ for all $k\in\N$, from (\ref{14.1}) we obtain
  that $(u(\cdot,t_k))_{k\in\N}$ forms a Cauchy sequence in $\big(W_N^{2,\infty}(\Om)\big)^\star$, and that hence (\ref{15.1})
  holds with some $u_\infty \in \big(W_N^{2,\infty}(\Om)\big)^\star$.
  The characterization in (\ref{15.2}) thereupon results from a second application of (\ref{14.1}), this time to the particular
  sequence 
  $(0,t,2t,...)$, 
  $t>0$, which namely ensures the existence of $c_1(K)>0$ such that under the hypotheses stated above we have
  \bas
	\|u(\cdot,t)-u_0\|_{(W_N^{2,\infty}(\Om))^\star} \le c_1(K) \io v_0
	\qquad \mbox{for all } t>0,
  \eas
  and thereby establishes (\ref{15.2}) due to (\ref{15.1}).
\qed
Also with regard to the large time behavior in the second solution component, we shall first content ourselves with a 
topological framework somewhat more moderate than the one appearing in Theorem \ref{theo18}:
\begin{lem}\label{lem166}
  If (\ref{init}) holds, then for $v$ as in Lemma \ref{lem13} we have
  \be{166.1}
	v(\cdot,t)\to 0
	\quad \mbox{in } L^1(\Omega)
	\qquad \mbox{as } t\to\infty.
  \ee
\end{lem}
\proof
  This can be seen by means of an argument similar to that performed to a slightly more complex variant in 
  \cite[Section 4]{win_TRAN2017}: 
  From Lemma \ref{lem3}, (\ref{vinfty}) and Lemma \ref{lem13} we obtain that $\int_0^\infty \io |\na v|^4$ is finite, and that hence,
  according to a Poincar\'e inequality,
  \bas
	\int_t^{t+1} \big\| v(\cdot,s)-\ov{v(\cdot,s)} \big\|_{L^4(\Om)} ds \to 0
	\qquad \mbox{as } t\to\infty,
  \eas
  where again $\ov{\vp}:=\frac{1}{|\Om|} \io \vp$ for $\vp\in L^1(\Om)$.
  Since furthermore $c_1:=\sup_{t>0} \|u(\cdot,t)\|_{L^\frac{4}{3}(\Om)}$ is finite by (\ref{17.1}), and since 
  \bas
	\int_t^{t+1} \io uv \to 0
	\qquad \mbox{as } t\to\infty
  \eas
  according to (\ref{13.67}), in view of the mass conservation property from (\ref{13.66}), and thanks to
  the H\"older inequality, this implies that
  \bas
	\ov{u_0} \int_t^{t+1} \|v(\cdot,s)\|_{L^1(\Om)} ds
	&=& \bigg| \int_t^{t+1} \io u(x,s) \ov{v(\cdot,s)} dxds \bigg| \\
	&=& \bigg| \int_t^{t+1} \io uv - \int_t^{t+1} \io u(x,s) \Big( v(x,s)-\ov{v(\cdot,s)} \Big) dxds \bigg| \\
	&\le& \int_t^{t+1} \io uv
	+ c_1 \int_t^{t+1} \big\| v(\cdot,s)-\ov{v(\cdot,s)} \big\|_{L^4(\Om)} ds \\[2mm]
	&\to& 0
	\qquad \mbox{as } t\to\infty.
  \eas
  Since $\ov{u_0}$ is positive according to (\ref{init}), from this we immediately infer (\ref{166.1}) upon 
  noting that $0\le t\mapsto \|v(\cdot,t)\|_{L^1(\Om)}$ is 
  nonincreasing due to the fact that $v$ solves its respective sub-problem in (\ref{0}) classically in $\Om\times (0,\infty)$
  by Lemma \ref{lem13}.
\qed
By means of straightforward interpolation relying on (\ref{18.1}) and (\ref{17.1}), from the latter and Lemma \ref{lem15}
we readily obtain our main result on stabilization in (\ref{0}):\abs
\proofc of Theorem \ref{theo18}. \quad
  Since (\ref{17.1}) guarantees that $(u(\cdot,t))_{t>0}$ is relatively compact with respect to the weak topology in each of the
  spaces $L^p(\Om)$ with $p>1$, taking $u_\infty$ as in Lemma \ref{lem15} we obtain the inclusion 
  $u_\infty \in \bigcap_{p\ge 1} L^p(\Om)$ and (\ref{18.1}) as direct consequences of (\ref{15.1}) when combined with the continuity
  of the embedding $(W^{2,\infty}_N(\Om))^\star \hra L^p(\Om)$ for any such $p$, while the identity 
  thereupon follows from (\ref{18.1}) and (\ref{13.66}).\\
  Likewise, (\ref{18.2}) results from the boundedness of $(v(\cdot,t))_{t>0}$ in $W^{1,\infty}(\Om)$, as implied by (\ref{17.1}),
  in conjunction with the statement on $L^1$ decay made in Lemma \ref{lem166}.
\qed
Returning to (\ref{15.2}), we can finally make sure that under a suitable smallness condition on $v_0$, the large time
limit thus obtained cannot be constant:
\begin{lem}\label{lem16}
  Let $u_0\in W^{1,\infty}(\Om)$ be nonnegative with $u_0\not\equiv const.$
  Then given $K>0$, one can find $\delta(K)>0$ such that if $v_0\in W^{1,\infty}(\Om)$ is nonnegative with 
  $\sqrt{v_0}\in W^{1,2}(\Om)$ and 
  $\|v_0\|_{L^\infty(\Om)} \le K$ as well as
  \be{16.1}
	\io v_0 \le \delta(K),
  \ee
  then the corresponding limit $u_\infty \in \big(W_N^{2,\infty}(\Om)\big)^\star$ from Lemma \ref{lem15} has the property that
  $u_\infty \not\equiv const.$
\end{lem}
\proof
  Since $u_0$ is continuous and not constant, we can fix numbers $c_1>0, c_2>c_1$ and $R>0$ as well as points $x_1\in\Om$ and 
  $x_2\in\Om$ such that $B_{2R}(x_i)\subset\Om$ for $i\in\{1,2\}$, and that $u_0\le c_1$ in $B_{2R}(x_1)$ and
  $u_0\ge c_2$ in $B_{2R}(x_2)$.
  It is then possible to pick $c_3>0$ as well as nonnegative functions $\psi_i\in C_0^\infty(\Om)$, $i\in\{1,2\}$, which are such that
  $\supp \psi_i \subset B_{2R}(x_i)$, that $\psi_i \equiv c_3$ in $B_R(x_i)$ and 
  $\|\psi_i\|_{W_N^{2,\infty}(\Om)} =1$
  for $i\in\{1,2\}$.
  For fixed $K>0$, we then take $\Xi(K)$ as in Lemma \ref{lem15} and claim that then the intended conclusion holds if we let
  \be{16.3}
	\delta(K):=\frac{c_3 \kappa \cdot |B_R(0)|}{2\Xi(K)}
  \ee
  with $\kappa:=\frac{c_2-c_1}{2}$.\abs
  Indeed, assuming on the contrary 
  that $0\le v_0\in W^{1,\infty}(\Om)$ with $\sqrt{v_0}\in W^{1,2}(\Om)$ and $\|v_0\|_{L^\infty(\Om)} \le K$
  satisfied (\ref{16.1}) but had the property that for the associated limit $u_\infty$ we had $u_\infty\equiv a$ for some $a\in\R$,
  by definition of $\kappa$ we would either have $a\le c_2-\kappa$ or $a\ge c_1+\kappa$.
  In the latter of these cases, however, we could use the localization features of $\psi_1$ together with (\ref{16.3}) to estimate
{\allowdisplaybreaks
  \bas
	\|u_\infty-u_0\|_{(W_N^{2,\infty}(\Om))^\star} 
	&\ge& \io (a-u_0)\cdot\psi_1 \\
	&=& \int_{B_{2R}(x_1)} (a-u_0)\cdot\psi_1 \\
	&\ge& \int_{B_R(x_1)} (a-u_0)\cdot c_3 \\
	&\ge& \int_{B_R(x_1)} \big\{ (c_1+\kappa)-c_1\big\}\cdot c_3 \\[2mm]
	&=& c_3\kappa \cdot |B_R(0)| \\[2mm]
	&=& 2\Xi(K) \delta(K),
  \eas
}
  which in view of (\ref{15.2}) is absurd.
  As it can be shown in quite a similar manner that also the inequality $a\le c_2-\kappa$ is impossible, it follows that, in fact,
  $u_\infty$ cannot coincide with any constant.
\qed
Our reasoning thereby becomes complete:\abs
\proofc of Theorem \ref{theo19}. \quad
  The claimed result has precisely been asserted by Lemma \ref{lem16}.
\qed

\bigskip

{\bf Acknowledgement.} \quad
  The author warmly thanks Wenbin Lv 
  as well as the anonymous reviewers
  for several useful comments.
  He furthermore acknowledges support of the {\em Deutsche Forschungsgemeinschaft} (Project No.~462888149).

\end{document}